\documentclass[12pt,a4paper]{article}
\usepackage{amsmath,amsthm,amssymb,latexsym,epic}

\newtheorem{theorem}{Theorem}
\newtheorem{lemma}[theorem]{Lemma}
\newtheorem{corollary}[theorem]{Corollary}
\newtheorem{proposition}[theorem]{Proposition}

\usepackage[all]{xy}
\usepackage{enumerate}
\numberwithin{equation}{section}

\newcommand{\N}{\mathbb{N}}

\renewcommand{\S}{\mathcal{S}}

\newcommand{\B}{\mathfrak{B}}

\renewcommand{\a}{{\alpha}}
\renewcommand{\b}{{\beta}}

\newcommand{\rank}{\mathrm{rank}}

\renewcommand{\phi}{\varphi}

\begin{document}

\title{$\mathcal{L}$-- and $\mathcal{R}$--cross-sections in the Brauer semigroup}
\author{Ganna Kudryavtseva, Victor Maltcev \\ and Volodymyr Mazorchuk}
\date{}
\maketitle

\begin{abstract}
We classify all cross-sections of Green's relations $\mathcal{L}$ and 
$\mathcal{R}$ in the Brauer semigroup. The regular behavior of such
cross-sections starts from $n=7$. We show that in the regular case
there are essentially two different cross-sections and all others are
$\S_n$-conjugated to one of these two. We also classify all 
cross-sections up to isomorphism. 
\end{abstract}

\section{Introduction}\label{s1}

Let $S$ be a semigroup, and $\rho$ be an equivalence relation on $S$.
A subsemigroup, $T\subset S$, is called a {\em cross-section} with respect
to $\rho$ provided that $T$ contains exactly one element from every 
equivalence class. Certainly the most natural equivalence relations on
a semigroup are congruences and Green's relations and so finding descriptions of
the cross-sections on these relations is a natural problem to consider. 

In what follows we will call the cross-sections with respect to Green's 
relations $\mathcal{L}$ ($\mathcal{R}$, $\mathcal{H}$, $\mathcal{D}$, 
$\mathcal{J}$) the $\mathcal{L}$-- ($\mathcal{R}$--, 
$\mathcal{H}$--, $\mathcal{D}$--, $\mathcal{J}$--) {\em cross-sections},
respectively. During the last decade the cross-sections of Green's relations
for some classical semigroups were studied by different authors. In particular,
for the symmetric inverse semigroup $\mathcal{IS}_n$ all $\mathcal{H}$--cross-sections 
were classified in \cite{CR} and all $\mathcal{L}$-- and $\mathcal{R}$--cross-sections 
were classified in \cite{GM3}. For the infinite symmetric inverse semigroup 
$\mathcal{IS}_X$ all $\mathcal{H}$--, $\mathcal{L}$-- and $\mathcal{R}$--cross-sections  
were classified in \cite{Pe2}, and for the full transformation semigroup 
$\mathcal{T}_n$ all $\mathcal{H}$-- and $\mathcal{R}$--cross-sections were classified 
in \cite{Pe1}. In \cite{Ml} it was shown that the description of 
$\mathcal{L}$-- and $\mathcal{R}$--cross-sections for the partial Brauer semigroup
$\mathcal{P}\mathfrak{B}\sb n$ reduces to the corresponding description for 
$\mathcal{IS}_n$, and a classification of $\mathcal{L}$-- and $\mathcal{R}$--cross-sections 
for the composition semigroup $\mathfrak{C}\sb n$ (defined in \cite{Mr} and studied in 
\cite{M}) was obtained.

The semigroups $\mathcal{P}\mathfrak{B}\sb n$ and $\mathfrak{C}\sb n$ are 
generalizations of the Brauer semigroup $\B_n$, introduced in \cite{Br}. In the 
present paper we give a classification of all 
$\mathcal{L}$-- and  $\mathcal{R}$--cross-section in $\B_n$. Surprisingly enough, for 
$\B_n$ the problem happens to be much more difficult than for $\mathcal{P}\mathfrak{B}\sb n$ 
and $\mathfrak{C}\sb n$. Because of the existence of an anti-involution on $\B_n$ it is 
enough to classify only one kind of cross-section. We do this for 
$\mathcal{R}$--cross-sections. It happens that $\mathcal{R}$--cross-sections exhibit a 
regular behavior starting with $n=7$. The symmetric group $\S_n\subset\B_n$ acts on the set 
of all $\mathcal{R}$--cross-sections of $\B_n$ in a natural way. For $n\geq 7$ we show that 
this action has exactly $2$ orbits, each containing $n!/2$ elements. We also describe the 
canonical representatives in these orbits, which we call the regular and the alternating
$\mathcal{R}$--cross-sections, respectively. We show that these two 
$\mathcal{R}$--cross-sections are not isomorphic as monoids. The cases $n\leq 6$ are
considered separately as the descriptions in these cases do not fit into the 
``regular'' picture.

The paper is organized as follows. In Section~\ref{s2} we give all the necessary 
background about the Brauer semigroups required for the sequel. In Section~\ref{s3} 
we study and completely determine a class of $\mathcal{R}$--cross-sections, which we call 
canonical. We use these results in Section~\ref{s4} to give a classification of all 
$\mathcal{R}$--cross-sections of $\B_n$. We also classify all $\mathcal{R}$--cross-sections
of $\B_n$ up to isomorphism. We finish the paper with a discussion of the
problems to classify $\mathcal{D}$-- and $\mathcal{H}$--cross-sections of $\B_n$
in Section~\ref{s5}.

\section{Preliminaries about the Brauer semigroup}\label{s2}

Let $n\in \N$, $M=\{1,2,\dots,n\}$ and $M'=\{1',2',\dots,n'\}$. We 
consider the map ${}':M\to M'$ as a fixed bijection and will denote the inverse 
bijection by the same symbol, that is $(a')'=a$.

Denote by $\B_n=\B(M)$ the set of all possible partitions of $M\cup M'$ into 
two-element subsets. It is a simple exercise to verify that $|\B_n|=(2n-1)!!$, 
see for example \cite{K}. For $\a\in \B_n$ and $a\neq b\in M\cup M'$ we set 
$a\equiv_{\a}b$ provided that $\{a,b\}\in\a$. That is $\equiv_{\a}$ is the equivalence
relation corresponding to the partition $\a$. Let $\alpha=X_1\cup \dots\cup X_n$ and 
$\beta=Y_1\cup \dots\cup Y_n$  be two elements from $\B_n$. Define a new equivalence 
relation, $\equiv$, on $M\cup M'$ as follows:
\begin{itemize}
\item for $a,b\in M$ we have $a\equiv b$ if and only if $a\equiv_{\alpha} b$
or there is a sequence, $c_1,\cdots,c_{2s}$, $s\geq 1$, of elements in $M$, 
such that $a\equiv_{\alpha} c'_1$, $c_1\equiv_{\beta} c_2$, $c'_2\equiv_{\alpha} 
c'_3,\dots,$ $c_{2s-1}\equiv_{\beta} c_{2s}$, and $c'_{2s}\equiv_{\alpha} b$;
\item for $a,b\in M$ we have $a'\equiv b'$ if and only if 
$a'\equiv_{\beta} b'$ or there is a sequence, $c_1,\cdots,c_{2s}$, $s\geq 1$, of 
elements in $M$, such that  $a'\equiv_{\beta} c_1$, $c'_1\equiv_{\alpha} c'_2$, 
$c_2\equiv_{\beta} c_3,\dots,$ $c'_{2s-1}\equiv_{\alpha} c'_{2s}$, and $c_{2s}\equiv_{\beta} b'$;
\item for $a,b\in M$ we have $a\equiv b'$ if and only if 
$b'\equiv a$ if and only if there is a 
sequence, $c_1,\cdots$, $c_{2s-1}$, $s\geq 1$, of elements in $M$, such that  $a\equiv_{\alpha} c'_1$,
$c_1\equiv_{\beta} c_2$, $c'_2\equiv_{\alpha} c'_3,\dots,$ $c'_{2s-2}\equiv_{\alpha} c'_{2s-1}$,
and $c_{2s-1}\equiv_{\beta} b'$.
\end{itemize}
It is easy to see that $\equiv$ determines a partition of $M\cup M'$ into
two-element subsets and so belongs to $\B_n$. 

\begin{figure}
\special{em:linewidth 0.4pt}
\unitlength 0.80mm
\linethickness{0.4pt}
\begin{picture}(150.00,65.00)
\drawline(15.33,15.33)(15.33,65.00)
\drawline(15.33,65.00)(45.00,65.00)
\drawline(45.00,65.00)(45.00,15.00)
\drawline(45.00,15.00)(15.33,15.00)
\drawline(15.33,15.00)(15.33,16.67)
\drawline(60.00,15.00)(60.00,65.00)
\drawline(60.00,65.00)(90.00,65.00)
\drawline(90.00,65.00)(90.00,15.00)
\drawline(90.00,15.00)(60.00,15.00)
\drawline(115.00,15.00)(115.00,65.00)
\drawline(115.00,65.00)(145.00,65.00)
\drawline(145.00,65.00)(145.00,15.00)
\drawline(145.00,15.00)(115.00,15.00)
\drawline(9.67,60.00)(15.33,60.00)
\drawline(15.33,55.00)(10.00,55.00)
\drawline(10.00,50.00)(15.33,50.00)
\drawline(15.33,45.00)(10.33,45.00)
\drawline(10.33,40.00)(15.33,40.00)
\drawline(15.33,35.00)(10.00,35.00)
\drawline(9.67,30.00)(15.33,30.00)
\drawline(15.33,25.00)(10.00,25.00)
\drawline(10.00,25.00)(10.00,25.00)
\drawline(10.00,20.00)(15.33,20.00)
\drawline(45.00,60.00)(50.00,60.00)
\drawline(50.00,55.00)(45.00,55.00)
\drawline(45.00,50.00)(50.00,50.00)
\drawline(50.00,45.00)(45.00,45.00)
\drawline(45.00,40.00)(49.67,40.00)
\drawline(49.67,35.00)(45.00,35.00)
\drawline(50.00,30.00)(45.00,30.00)
\drawline(45.00,25.00)(50.00,25.00)
\drawline(50.00,20.00)(45.00,20.00)
\drawline(55.00,60.00)(60.00,60.00)
\drawline(60.00,55.00)(55.00,55.00)
\drawline(55.00,50.00)(60.00,50.00)
\drawline(60.00,45.00)(55.00,45.00)
\drawline(55.00,40.00)(60.00,40.00)
\drawline(60.00,35.00)(55.00,35.00)
\drawline(55.00,30.00)(60.00,30.00)
\drawline(60.00,25.00)(55.33,25.00)
\drawline(55.00,20.00)(60.00,20.00)
\drawline(90.00,60.00)(95.00,60.00)
\drawline(95.00,55.00)(90.00,55.00)
\drawline(90.00,50.00)(95.00,50.00)
\drawline(95.00,45.00)(90.00,45.00)
\drawline(90.00,40.00)(95.00,40.00)
\drawline(95.00,35.00)(90.00,35.00)
\drawline(90.00,30.00)(95.00,30.00)
\drawline(95.00,25.00)(90.00,25.00)
\drawline(90.00,20.00)(95.00,20.00)
\drawline(110.00,60.00)(115.00,60.00)
\drawline(115.00,55.00)(110.00,55.00)
\drawline(110.00,50.00)(115.00,50.00)
\drawline(115.00,45.00)(110.00,45.00)
\drawline(110.00,40.00)(115.00,40.00)
\drawline(115.00,35.00)(110.00,35.00)
\drawline(110.00,30.00)(115.00,30.00)
\drawline(115.00,25.00)(110.00,25.00)
\drawline(110.00,20.00)(115.00,20.00)
\drawline(145.00,60.00)(150.00,60.00)
\drawline(150.00,55.00)(145.00,55.00)
\drawline(145.00,50.00)(150.00,50.00)
\drawline(150.00,45.00)(145.00,45.00)
\drawline(145.00,40.00)(150.00,40.00)
\drawline(150.00,35.00)(145.00,35.00)
\drawline(145.00,30.00)(150.00,30.00)
\drawline(150.00,25.00)(145.00,25.00)
\drawline(145.00,20.00)(150.00,20.00)
\drawline(15.33,20.00)(45.00,40.00)
\drawline(45.00,35.00)(42.33,35.00)
\drawline(42.33,35.00)(42.33,30.00)
\drawline(42.33,30.00)(45.00,30.00)
\drawline(45.00,25.00)(42.33,25.00)
\drawline(42.33,25.00)(42.33,20.00)
\drawline(42.33,20.00)(45.00,20.00)
\drawline(45.00,50.00)(42.33,50.00)
\drawline(42.33,50.00)(42.33,45.00)
\drawline(42.33,45.00)(45.00,45.00)
\drawline(60.00,40.00)(63.00,40.00)
\drawline(63.00,40.00)(63.00,35.00)
\drawline(63.00,35.00)(60.00,35.00)
\drawline(64.67,45.00)(60.00,45.00)
\drawline(60.00,50.00)(90.00,30.00)
\drawline(15.33,45.00)(45.00,60.00)
\drawline(15.33,55.00)(45.00,55.00)
\drawline(15.33,40.00)(18.33,40.00)
\drawline(18.33,40.00)(18.33,30.00)
\drawline(18.33,30.00)(15.33,30.00)
\drawline(15.33,25.00)(20.33,25.00)
\drawline(20.33,25.00)(20.33,35.00)
\drawline(20.33,35.00)(15.33,35.00)
\drawline(15.33,50.00)(18.33,50.00)
\drawline(18.33,50.00)(18.33,60.00)
\drawline(18.33,60.00)(15.33,60.00)
\drawline(60.00,60.00)(63.00,60.00)
\drawline(63.00,60.00)(63.00,55.00)
\drawline(63.00,55.00)(60.00,55.00)
\drawline(90.00,60.00)(87.00,60.00)
\drawline(87.00,60.00)(87.00,35.00)
\drawline(87.00,35.00)(90.00,35.00)
\drawline(90.00,25.00)(87.00,25.00)
\drawline(87.00,25.00)(87.00,20.00)
\drawline(87.00,20.00)(90.00,20.00)
\drawline(90.00,40.00)(88.33,40.00)
\drawline(88.33,40.00)(88.33,45.00)
\drawline(88.33,45.00)(90.00,45.00)
\drawline(90.00,50.00)(88.33,50.00)
\drawline(88.33,50.00)(88.33,55.00)
\drawline(88.33,55.00)(90.00,55.00)
\drawline(115.00,60.00)(119.00,60.00)
\drawline(119.00,60.00)(119.00,50.00)
\drawline(119.00,50.00)(115.00,50.00)
\drawline(115.00,55.00)(117.67,55.00)
\drawline(117.67,55.00)(117.67,45.00)
\drawline(117.67,45.00)(115.00,45.00)
\drawline(115.00,40.00)(117.67,40.00)
\drawline(117.67,40.00)(117.67,30.00)
\drawline(117.67,30.00)(115.00,30.00)
\drawline(115.00,35.00)(119.00,35.00)
\drawline(119.00,35.00)(119.00,25.00)
\drawline(119.00,25.00)(115.00,25.00)
\drawline(115.00,20.00)(145.00,30.00)
\drawline(145.00,25.00)(141.33,25.00)
\drawline(141.33,25.00)(141.33,20.00)
\drawline(141.33,20.00)(145.00,20.00)
\drawline(145.00,60.00)(141.33,60.00)
\drawline(141.33,60.00)(141.33,35.00)
\drawline(141.33,35.00)(145.00,35.00)
\drawline(145.00,40.00)(143.00,40.00)
\drawline(143.00,40.00)(143.00,45.00)
\drawline(143.00,45.00)(145.00,45.00)
\drawline(145.00,50.00)(143.00,50.00)
\drawline(143.00,50.00)(143.00,55.00)
\drawline(143.00,55.00)(145.00,55.00)
\put(102.33,40.00){\makebox(0,0)[cc]{$=$}}
\put(52.33,40.00){\makebox(0,0)[cc]{$\cdot$}}
\drawline(65.00,45.00)(65.00,30.00)
\drawline(60.00,25.00)(63.00,25.00)
\drawline(63.00,25.00)(63.00,20.00)
\drawline(63.00,20.00)(60.00,20.00)
\drawline(65.00,30.00)(60.00,30.00)
\end{picture}
\caption{Chips and their multiplication.}\label{fig:f1}
\end{figure}
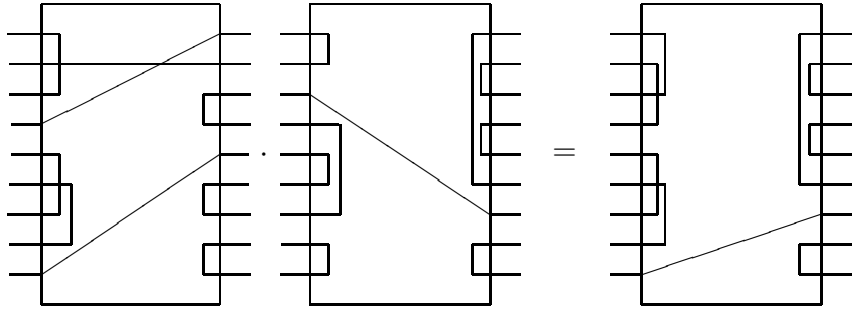

One can think about the elements from $\B_n$ as certain ``microchips'' with $n$ pins on the
left hand side (corresponding to $M$) and $n$ pins on the right hand side (corresponding to $M'$).
Having  $\alpha\in \B_n$ we connect two pins in the corresponding chip if and only if they belong
to the same set of the partition $\alpha$. The operation described above can then be viewed as
a ``composition'' of such chips: having $\alpha,\beta\in\B_n$ we identify (connect) the right pins 
of $\alpha$ with the corresponding left pins of $\beta$, which uniquely defines a connection of
the remaining pins (which are the left pins of $\alpha$ and the right pins of $\beta$). An example of multiplication of two chips from $\B_n$ is given on Figure~\ref{fig:f1}. Note that 
performing the operation we can obtain some ``dead circles'' formed by some identified pins from
$\alpha$ and $\beta$, see for example the two lowest identified pins on Figure~\ref{fig:f1}. 
These circles should be disregarded. From such an interpretation it is not hard to see 
that the composition of elements from $\B_n$ defined above is associative (see \cite{Br}).
 
The obtained semigroup is called the {\em Brauer} semigroup (monoid). The (deformed) semigroup 
algebra of $\B_n$ is the famous Brauer algebra, which was introduced in \cite{Br}, and which 
plays an important role in the study of certain representations of orthogonal groups. The
Brauer algebra has been extensively studied in the literature (see for example
\cite{KX} and the references therein). The Brauer semigroup was studied in, for example,
\cite{K,M2,M,Ml2}.

There is a natural monomorphism $\S_n\hookrightarrow \B_n$ defined as follows:
\begin{displaymath}
\sigma\mapsto \{1,\sigma(1)'\}\cup \{2,\sigma(2)'\}\cup\dots\cup
\{n,\sigma(n)'\},\quad\text{where}\quad \sigma\in \S_n,
\end{displaymath}
and we will identify $\sigma\in\S_n$ with its image in $\B_n$ via this embedding.
The image of $\S_n$ under the above embedding coincides with the set of all
invertible elements in $\B_n$. For $n\geq 2k$ there is also a natural monomorphism
from the symmetric inverse semigroup $\mathcal{IS}_k$ on $\{1,\dots,k\}$ into
$\B_n$, which is  constructed in the following way: an element, 
$\sigma\in \mathcal{IS}_k$, is sent to the element $\gamma\in \B_n$, uniquely defined
by the following conditions: 
\begin{itemize}
\item $i\equiv_{\gamma}i'$ for all $i>2k$;
\item $2i\equiv_{\gamma}(2\sigma(i))'$ and $2i-1\equiv_{\gamma}(2\sigma(i)-1)'$
for all $i\in \{1,\dots,k\}$, which belong to the domain of $\sigma$;
\item $2i\equiv_{\gamma}2i-1$ for all $i\in \{1,\dots,k\}$, which do not belong 
to the domain of $\sigma$;
\item $(2i)'\equiv_{\gamma}(2i-1)'$ for all $i\in \{1,\dots,k\}$, which do not belong 
to the range of $\sigma$.
\end{itemize}

For $\a\in\B_n$ a partition set, $X\in\a$, having the form $X=\{a,b'\}$ for $a,b\in M$,
will be called a {\em line} in $\a$. The number of different lines in $\a$
is called the {\em rank} of $\a$ and is denoted by $\mathrm{rank}(\a)$. The
number $n-\mathrm{rank}(\a)$ is called the {\em corank} of $\a$ and is denoted by 
$\mathrm{corank}(\a)$. For example, $\sigma\in\S_n$ if and only if 
$\mathrm{rank}(\sigma)=n$ if and only if $\mathrm{corank}(\sigma)=0$. Note that 
$\mathrm{corank}(\a)$ is even for every $\a\in\B_n$ and that
$\mathrm{corank}(\a\beta)\leq \mathrm{corank}(\a)+\mathrm{corank}(\beta)$
for all $\a,\beta\in\B_n$. Since $\B_n$ is finite, for 
every $\a\in\B_n$ there exists a unique idempotent, $\pi\in\B_n$, such that
$\pi=\a^i$ for some $i\in\N$. The rank of $\pi$ is thus an invariant for 
$\a$ and will be called the {\em stable rank} $\mathrm{strank}(\a)$ of $\a$. Note
that $\mathrm{strank}(\a)\leq \mathrm{rank}(\a)$.

Call $\a,\b\in\B_n$ {\em left neighbors} or {\em right neighbors} provided
that $\{a,b\}\in\a$ if and only if $\{a,b\}\in\b$ for all $a,b\in M$ and
for all $a,b\in M'$, respectively. Green's relations $\mathcal{L}$, 
$\mathcal{R}$, $\mathcal{D}$, $\mathcal{H}$, $\mathcal{J}$, in $\B_n$ can be
described as follows:

\begin{theorem}\label{tgreen}\cite[Theorem~7]{M2}
Let $\a,\b\in \B_n$. Then
\begin{enumerate}[(i)]
\item\label{tgreen.1} $\a\mathcal{L}\b$ if and only if $\a$ and $\b$ are right
neighbors if and only if there exists $\sigma\in\S_n$ such that $\sigma\a=\b$;
\item\label{tgreen.2} $\a\mathcal{R}\b$ if and only if $\a$ and $\b$ are left
neighbors if and only if there exists $\sigma\in\S_n$ such that $\a\sigma=\b$;
\item\label{tgreen.3} $\a\mathcal{H}\b$ if and only if $\a$ and $\b$ are both
left and right neighbors;
\item\label{tgreen.4} $\a\mathcal{D}\b$ if and only if $\mathrm{rank}(\a)=
\mathrm{rank}(\b)$ if and only if there exist $\sigma,\tau\in\S_n$ such that 
$\sigma\a\tau=\b$;
\item\label{tgreen.5} $\mathcal{D}=\mathcal{J}$.
\end{enumerate}
\end{theorem}

A subset, $X\subset M$, will be called {\em $\alpha$-invariant} for some $\alpha\in\B_n$ 
provided that for any $a\in X\cup X'$ and any $b\in M\cup M'$ the
condition $a\equiv_{\alpha} b$ implies  $b\in X\cup X'$. If $X$ is invariant
with respect to $\alpha$, then we define the element $\alpha|_{X}\in\B(X)$ in the
following way: 
\begin{itemize}
\item for all $a,b\in X\cup X'$ we have $a\equiv_{\alpha|_{X}} b$ if and
only if $a\equiv_{\alpha} b$.
\end{itemize}
The element $\alpha|_{X}$ is called the {\em restriction} of $\alpha$ to $X$. 
Note that if $X$ is $\alpha$-invariant then $M\setminus X$ is $\alpha$-invariant 
as well.

The involution ${}':M\cup M'\to M\cup M'$ extends in a natural way to the
anti-involution ${}^*:\B_n\to \B_n$, which, in the language of chips, acts on
a chip by taking the mirror image of it. It is obvious that 
$\a\a^*\a=\a$, in particular, $\B_n$ is a regular semigroup. Since ${}^*$ is
an anti-involution, it interchanges $\mathcal{L}$-- and $\mathcal{R}$--classes
and the corresponding cross-sections. Hence, it is enough to classify one type
of cross-sections, say the $\mathcal{R}$--cross-sections. The classification of
the $\mathcal{L}$--cross-sections is then obtained by applying ${}^*$.

\section{Canonical $\mathcal{R}$--cross-sections in $\B_n$}\label{s3}

It does not follow from the definition that $\mathcal{R}$--cross-sections in $\B_n$
exist. In this section we construct and investigate a special (rather big) family of 
$\mathcal{R}$--cross-sections in $\B_n$, in particular, showing that they exist.
We call an $\mathcal{R}$--cross-section, $\Lambda$,  of $\B_n$ {\em canonical}
provided that for every  $\a\in\Lambda$ there exists $k\in\{0,1,\dots,
\lfloor\frac{n}{2}\rfloor\}$ such that
$\mathrm{corank}(\a)=2k$ and 
\begin{displaymath}
(n-2k+1)'\equiv_{\a} (n-2k+2)', (n-2k+3)'\equiv_{\a} 
(n-2k+4)',\dots,  (n-1)'\equiv_{\a} n'.
\end{displaymath}
Our aim in this section is to show that
canonical $\mathcal{R}$--cross-sections of $\B_n$ exist and classify 
all such cross-sections. However, even the existence is not obvious and will be
established only in Proposition~\ref{p10} and Proposition~\ref{p12}. To be able 
to prove these results in the early part of the section we describe what must happen 
within a canonical $\mathcal{R}$--cross-sections of $\B_n$, should one exist.

Let $\Lambda$ be a canonical $\mathcal{R}$--cross-section in $\B_n$. For 
$1\leq i< j\leq n$ denote by $\alpha_{i,j}$ the (unique) element of $\Lambda$
such that $\mathrm{corank}(\alpha_{i,j})=2$ and $i\equiv_{\alpha_{i,j}} j$.
For $k=1,\dots,\lfloor\frac{n}{2}\rfloor$ set
\begin{displaymath}
\Lambda_k=\{\a\in\Lambda\,:\, \mathrm{corank}(\a)=2k\}.
\end{displaymath}
Note that $\Lambda_0=\{\mathrm{id}\}$ and 
$\Lambda_1=\{\a_{i,j}\,:\, 1\leq i< j\leq n\}$. Hence for
$k=1,\dots,\lfloor\frac{n}{2}\rfloor$ we can define
\begin{displaymath}
\Lambda_1^{(k)}=\{\a_{i,j}\,:\, j\leq n-2(k-1)\}.
\end{displaymath}
In particular, $\Lambda_1^{(1)}=\Lambda_1$.
Later on we will show that the elements $\alpha_{i,j}$ completely determine 
(generate) $\Lambda$. Our main idea is to collect enough information 
(mainly technical) about $\alpha_{i,j}$ to be able to explicitly describe $\Lambda$. 

\begin{figure}
\special{em:linewidth 0.4pt}
\unitlength 0.80mm
\linethickness{0.4pt}
\begin{picture}(150.00,65.00)
\drawline(115.00,15.00)(115.00,65.00)
\drawline(115.00,65.00)(145.00,65.00)
\drawline(145.00,65.00)(145.00,15.00)
\drawline(145.00,15.00)(115.00,15.00)
\drawline(110.00,60.00)(115.00,60.00)
\drawline(115.00,55.00)(110.00,55.00)
\drawline(110.00,50.00)(115.00,50.00)
\drawline(115.00,45.00)(110.00,45.00)
\drawline(110.00,40.00)(115.00,40.00)
\drawline(115.00,35.00)(110.00,35.00)
\drawline(110.00,30.00)(115.00,30.00)
\drawline(115.00,25.00)(110.00,25.00)
\drawline(110.00,20.00)(115.00,20.00)
\drawline(145.00,60.00)(150.00,60.00)
\drawline(150.00,55.00)(145.00,55.00)
\drawline(145.00,50.00)(150.00,50.00)
\drawline(150.00,45.00)(145.00,45.00)
\drawline(145.00,40.00)(150.00,40.00)
\drawline(150.00,35.00)(145.00,35.00)
\drawline(145.00,30.00)(150.00,30.00)
\drawline(150.00,25.00)(145.00,25.00)
\drawline(145.00,20.00)(150.00,20.00)
\drawline(115.00,20.00)(145.00,20.00)
\drawline(115.00,30.00)(145.00,30.00)
\drawline(115.00,35.00)(145.00,35.00)
\drawline(115.00,40.00)(145.00,40.00)
\drawline(115.00,45.00)(145.00,45.00)
\drawline(115.00,50.00)(145.00,50.00)
\drawline(145.00,60.00)(140.00,60.00)
\drawline(145.00,55.00)(140.00,55.00)
\drawline(140.00,55.00)(140.00,60.00)
\drawline(115.00,25.00)(122.00,25.00)
\drawline(115.00,60.00)(122.00,60.00)
\drawline(122.00,25.00)(122.00,60.00)
\drawline(115.00,55.00)(145.00,25.00)
\drawline(65.00,15.00)(65.00,65.00)
\drawline(65.00,65.00)(95.00,65.00)
\drawline(95.00,65.00)(95.00,15.00)
\drawline(95.00,15.00)(65.00,15.00)
\drawline(60.00,60.00)(65.00,60.00)
\drawline(65.00,55.00)(60.00,55.00)
\drawline(60.00,50.00)(65.00,50.00)
\drawline(65.00,45.00)(60.00,45.00)
\drawline(60.00,40.00)(65.00,40.00)
\drawline(65.00,35.00)(60.00,35.00)
\drawline(60.00,30.00)(65.00,30.00)
\drawline(65.00,25.00)(60.00,25.00)
\drawline(60.00,20.00)(65.00,20.00)
\drawline(95.00,60.00)(100.00,60.00)
\drawline(100.00,55.00)(95.00,55.00)
\drawline(95.00,50.00)(100.00,50.00)
\drawline(100.00,45.00)(95.00,45.00)
\drawline(95.00,40.00)(100.00,40.00)
\drawline(100.00,35.00)(95.00,35.00)
\drawline(95.00,30.00)(100.00,30.00)
\drawline(100.00,25.00)(95.00,25.00)
\drawline(95.00,20.00)(100.00,20.00)
\drawline(65.00,25.00)(95.00,25.00)
\drawline(65.00,30.00)(95.00,30.00)
\drawline(65.00,20.00)(95.00,20.00)
\drawline(65.00,40.00)(95.00,40.00)
\drawline(65.00,45.00)(95.00,45.00)
\drawline(65.00,50.00)(95.00,50.00)
\drawline(95.00,55.00)(90.00,55.00)
\drawline(95.00,60.00)(90.00,60.00)
\drawline(90.00,55.00)(90.00,60.00)
\drawline(65.00,35.00)(67.00,35.00)
\drawline(65.00,55.00)(67.00,55.00)
\drawline(67.00,55.00)(67.00,35.00)
\drawline(65.00,60.00)(95.00,35.00)
\drawline(15.00,15.00)(15.00,65.00)
\drawline(15.00,65.00)(45.00,65.00)
\drawline(45.00,65.00)(45.00,15.00)
\drawline(45.00,15.00)(15.00,15.00)
\drawline(10.00,60.00)(15.00,60.00)
\drawline(15.00,55.00)(10.00,55.00)
\drawline(10.00,50.00)(15.00,50.00)
\drawline(15.00,45.00)(10.00,45.00)
\drawline(10.00,40.00)(15.00,40.00)
\drawline(15.00,35.00)(10.00,35.00)
\drawline(10.00,30.00)(15.00,30.00)
\drawline(15.00,25.00)(10.00,25.00)
\drawline(10.00,20.00)(15.00,20.00)
\drawline(45.00,60.00)(50.00,60.00)
\drawline(50.00,55.00)(45.00,55.00)
\drawline(45.00,50.00)(50.00,50.00)
\drawline(50.00,45.00)(45.00,45.00)
\drawline(45.00,40.00)(50.00,40.00)
\drawline(50.00,35.00)(45.00,35.00)
\drawline(45.00,30.00)(50.00,30.00)
\drawline(50.00,25.00)(45.00,25.00)
\drawline(45.00,20.00)(50.00,20.00)
\drawline(15.00,20.00)(45.00,20.00)
\drawline(15.00,25.00)(45.00,25.00)
\drawline(15.00,30.00)(45.00,30.00)
\drawline(15.00,35.00)(20.00,35.00)
\drawline(15.00,40.00)(20.00,40.00)
\drawline(20.00,35.00)(20.00,40.00)
\drawline(15.00,45.00)(20.00,45.00)
\drawline(15.00,50.00)(20.00,50.00)
\drawline(20.00,45.00)(20.00,50.00)
\drawline(15.00,55.00)(20.00,55.00)
\drawline(15.00,60.00)(20.00,60.00)
\drawline(20.00,55.00)(20.00,60.00)
\drawline(45.00,55.00)(40.00,55.00)
\drawline(45.00,60.00)(40.00,60.00)
\drawline(40.00,55.00)(40.00,60.00)
\drawline(45.00,35.00)(40.00,35.00)
\drawline(45.00,40.00)(40.00,40.00)
\drawline(40.00,35.00)(40.00,40.00)
\drawline(45.00,45.00)(40.00,45.00)
\drawline(45.00,50.00)(40.00,50.00)
\drawline(40.00,45.00)(40.00,50.00)
\put(30.00,10.00){\makebox(0,0)[cc]{$\beta_{3}$}}
\put(80.00,10.00){\makebox(0,0)[cc]{$\alpha_{4,8}$}}
\put(130.00,10.00){\makebox(0,0)[cc]{$\alpha_{2,9}$}}
\end{picture}
\caption{Some idempotents in $\Lambda$ for $\B_9$.}\label{fig:f2}
\end{figure}
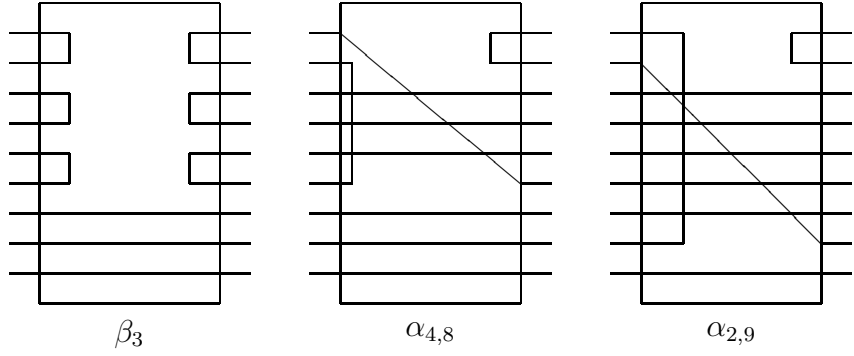

For every $i=0,\dots,\lfloor\frac{n}{2}\rfloor$ define the element 
$\beta_{i}$ of corank $2i$ as follows: $j\equiv_{\beta_{i}}j'$ for all 
$j=1,\dots,n-2i$; $j\equiv_{\beta_{i}}j+1$ and $j'\equiv_{\beta_{i}}(j+1)'$ 
for all $j=n-2i+1,n-2i+3,\dots,n-1$ (see example on Figure~\ref{fig:f2}).
We start with the following observation:

\begin{lemma}\label{l1}
\begin{enumerate}[(i)]
\item\label{l1.1} Let $\a\in\Lambda$ and suppose that there exists $\pi\in\B_n$ 
such that $\a\mathcal{H}\pi$ and $\pi^2=\pi$. Then $\a=\pi$.
\item\label{l1.n2}  For every $i=0,\dots,\lfloor\frac{n}{2}\rfloor$,
the element $\beta_{i}$ is an idempotent, and so belongs to $\Lambda$.
\item\label{l1.n3} For $j=n-1,n$ and for all $i=1,\dots,j-1$ the element $\alpha_{i,j}$
of $\Lambda$ is an idempotent and satisfies $s\equiv_{\alpha_{i,j}}s'$ for all 
$s\neq i,n-1,n$. Moreover,
in the case $(i,j)\neq(n-1,n)$ we also have $i'\equiv_{\alpha_{i,j}}\overline{j}$, 
where $\{j,\overline{j}\}=\{n-1,n\}$ (see example on Figure~\ref{fig:f2}). In the
case $(i,j)=(n-1,n)$ we have $\alpha_{i,j}=\beta_1$.
\item\label{l1.2} The element $\alpha_{i,j}$ is an idempotent if and only if $j=n-1$
or $j=n$.
\item\label{l1.n5} If two elements of $\Lambda$ are $\mathcal{D}$-related then they
are $\mathcal{L}$-related.
\end{enumerate}
\end{lemma}

\begin{proof}
Since $\pi$ is an idempotent, its $\mathcal{H}$-class is a (finite) maximal subgroup 
of $\B_n$. In particular, $\a^i=\pi$ for some $i\in\N$. This implies $\pi\in\Lambda$.
Furthermore, $\a\mathcal{H}\pi$ implies $\a\mathcal{R}\pi$ and hence $\pi=\a$ as
$\Lambda$ is an $\mathcal{R}$--cross-section. This proves \eqref{l1.1}. 

For any $i$, $\beta_{i}$ is easily seen to be an idempotent by direct calculation.
If $\alpha\in\Lambda$ and $\alpha\mathcal{R}\beta_{i}$ then we have 
$\alpha\mathcal{H}\beta_{i}$ by Theorem~\ref{tgreen}\eqref{tgreen.3}, and so
$\alpha=\beta_i$ by part \eqref{l1.1}. This proves \eqref{l1.n2}. 
\eqref{l1.n3} and \eqref{l1.2} are proved by direct calculation, and
\eqref{l1.n5} follows from the definition of a canonical
$\mathcal{R}$--cross-section anf Theorem~\ref{tgreen}\eqref{tgreen.4}.
\end{proof}

\begin{corollary}\label{c3}
Let $1\leq i<j\leq n-2$.
\begin{enumerate}[(i)]
\item\label{c3.1} Assume that $j$ and $n$ have the same parity. Then
$s\equiv_{\alpha_{i,j}}s'$ for all $s\in\{1,\dots,j-2\}\setminus\{i\}$,
and, if $i\neq j-1$, we also have $(j-1)\equiv_{\alpha_{i,j}}i'$.
\item\label{c3.2} Assume that $j$ and $n$ have different parities. Then
$s\equiv_{\alpha_{i,j}}s'$ for all $s\in\{1,\dots,j-1\}\setminus\{i\}$,
and $(j+1)\equiv_{\alpha_{i,j}}i'$.
\item\label{c3.3} For every $k=1,\dots,\lfloor\frac{n-j}{2}\rfloor$
and for every $s\in\{n-2k+2,n-2k+1\}$ there exists $t\in \{n-2k,n-2k-1\}$
such that $s\equiv_{\alpha_{i,j}}t'$.
\end{enumerate}
\end{corollary}

\begin{figure}
\special{em:linewidth 0.4pt}
\unitlength 0.80mm
\linethickness{0.4pt}
\begin{picture}(150.00,65.00)
\drawline(15.00,15.00)(15.00,65.00)
\drawline(15.00,65.00)(45.00,65.00)
\drawline(45.00,65.00)(45.00,15.00)
\drawline(45.00,15.00)(15.00,15.00)
\drawline(15.00,20.00)(10.00,20.00)
\drawline(15.00,25.00)(10.00,25.00)
\drawline(15.00,30.00)(10.00,30.00)
\drawline(15.00,35.00)(10.00,35.00)
\drawline(15.00,40.00)(10.00,40.00)
\drawline(15.00,45.00)(10.00,45.00)
\drawline(15.00,50.00)(10.00,50.00)
\drawline(15.00,55.00)(10.00,55.00)
\drawline(15.00,60.00)(10.00,60.00)
\drawline(45.00,20.00)(50.00,20.00)
\drawline(45.00,25.00)(50.00,25.00)
\drawline(45.00,30.00)(50.00,30.00)
\drawline(45.00,35.00)(50.00,35.00)
\drawline(45.00,40.00)(50.00,40.00)
\drawline(45.00,45.00)(50.00,45.00)
\drawline(45.00,50.00)(50.00,50.00)
\drawline(45.00,55.00)(50.00,55.00)
\drawline(45.00,60.00)(50.00,60.00)
\drawline(45.00,20.00)(15.00,20.00)
\drawline(45.00,25.00)(15.00,25.00)
\drawline(45.00,30.00)(15.00,30.00)
\drawline(15.00,35.00)(20.00,35.00)
\drawline(15.00,40.00)(20.00,40.00)
\drawline(20.00,35.00)(20.00,40.00)
\drawline(15.00,45.00)(20.00,45.00)
\drawline(15.00,50.00)(20.00,50.00)
\drawline(20.00,45.00)(20.00,50.00)
\drawline(15.00,55.00)(20.00,55.00)
\drawline(15.00,60.00)(20.00,60.00)
\drawline(20.00,55.00)(20.00,60.00)
\drawline(45.00,35.00)(40.00,35.00)
\drawline(45.00,40.00)(40.00,40.00)
\drawline(40.00,35.00)(40.00,40.00)
\drawline(45.00,45.00)(40.00,45.00)
\drawline(45.00,50.00)(40.00,50.00)
\drawline(40.00,45.00)(40.00,50.00)
\drawline(45.00,55.00)(40.00,55.00)
\drawline(45.00,60.00)(40.00,60.00)
\drawline(40.00,55.00)(40.00,60.00)
\drawline(60.00,15.00)(60.00,65.00)
\drawline(60.00,65.00)(90.00,65.00)
\drawline(90.00,65.00)(90.00,15.00)
\drawline(90.00,15.00)(60.00,15.00)
\drawline(55.00,20.00)(60.00,20.00)
\drawline(55.00,25.00)(60.00,25.00)
\drawline(55.00,30.00)(60.00,30.00)
\drawline(55.00,35.00)(60.00,35.00)
\drawline(55.00,40.00)(60.00,40.00)
\drawline(55.00,45.00)(60.00,45.00)
\drawline(55.00,50.00)(60.00,50.00)
\drawline(55.00,55.00)(60.00,55.00)
\drawline(55.00,60.00)(60.00,60.00)
\drawline(95.00,20.00)(90.00,20.00)
\drawline(95.00,25.00)(90.00,25.00)
\drawline(95.00,30.00)(90.00,30.00)
\drawline(95.00,35.00)(90.00,35.00)
\drawline(95.00,40.00)(90.00,40.00)
\drawline(95.00,45.00)(90.00,45.00)
\drawline(95.00,50.00)(90.00,50.00)
\drawline(95.00,55.00)(90.00,55.00)
\drawline(95.00,60.00)(90.00,60.00)
\drawline(60.00,20.00)(90.00,20.00)
\drawline(60.00,30.00)(90.00,30.00)
\drawline(60.00,35.00)(90.00,25.00)
\drawline(60.00,25.00)(65.00,25.00)
\drawline(60.00,40.00)(65.00,40.00)
\drawline(65.00,25.00)(65.00,40.00)
\drawline(90.00,60.00)(85.00,60.00)
\drawline(90.00,55.00)(85.00,55.00)
\drawline(85.00,60.00)(85.00,55.00)
\dashline{1}(60.00,60.00)(90.00,50.00)
\dashline{1}(60.00,55.00)(90.00,45.00)
\dashline{2}(60.00,60.00)(90.00,45.00)
\dashline{2}(60.00,55.00)(90.00,50.00)
\dashline{3}(60.00,50.00)(90.00,40.00)
\dashline{3}(60.00,45.00)(90.00,35.00)
\dashline{4}(60.00,50.00)(90.00,35.00)
\dashline{5}(60.00,45.00)(90.00,40.00)
\drawline(115.00,15.00)(115.00,65.00)
\drawline(115.00,65.00)(145.00,65.00)
\drawline(145.00,65.00)(145.00,15.00)
\drawline(145.00,15.00)(115.00,15.00)
\drawline(115.00,20.00)(110.00,20.00)
\drawline(115.00,25.00)(110.00,25.00)
\drawline(115.00,30.00)(110.00,30.00)
\drawline(115.00,35.00)(110.00,35.00)
\drawline(115.00,40.00)(110.00,40.00)
\drawline(115.00,45.00)(110.00,45.00)
\drawline(115.00,50.00)(110.00,50.00)
\drawline(115.00,55.00)(110.00,55.00)
\drawline(115.00,60.00)(110.00,60.00)
\drawline(145.00,20.00)(150.00,20.00)
\drawline(145.00,25.00)(150.00,25.00)
\drawline(145.00,30.00)(150.00,30.00)
\drawline(145.00,35.00)(150.00,35.00)
\drawline(145.00,40.00)(150.00,40.00)
\drawline(145.00,45.00)(150.00,45.00)
\drawline(145.00,50.00)(150.00,50.00)
\drawline(145.00,55.00)(150.00,55.00)
\drawline(145.00,60.00)(150.00,60.00)
\drawline(145.00,20.00)(115.00,20.00)
\drawline(145.00,25.00)(115.00,25.00)
\drawline(145.00,30.00)(115.00,30.00)
\drawline(115.00,35.00)(120.00,35.00)
\drawline(115.00,40.00)(120.00,40.00)
\drawline(120.00,35.00)(120.00,40.00)
\drawline(115.00,45.00)(120.00,45.00)
\drawline(115.00,50.00)(120.00,50.00)
\drawline(120.00,45.00)(120.00,50.00)
\drawline(115.00,55.00)(120.00,55.00)
\drawline(115.00,60.00)(120.00,60.00)
\drawline(120.00,55.00)(120.00,60.00)
\drawline(145.00,35.00)(140.00,35.00)
\drawline(145.00,40.00)(140.00,40.00)
\drawline(140.00,35.00)(140.00,40.00)
\drawline(145.00,45.00)(140.00,45.00)
\drawline(145.00,50.00)(140.00,50.00)
\drawline(140.00,45.00)(140.00,50.00)
\drawline(145.00,55.00)(140.00,55.00)
\drawline(145.00,60.00)(140.00,60.00)
\drawline(140.00,55.00)(140.00,60.00)
\put(102.33,40.00){\makebox(0,0)[cc]{$=$}}
\put(52.33,40.00){\makebox(0,0)[cc]{$\cdot$}}
\put(30.00,10.00){\makebox(0,0)[cc]{$\beta_{3}$}}
\put(75.00,10.00){\makebox(0,0)[cc]{$\alpha_{2,5}$}}
\put(130.00,10.00){\makebox(0,0)[cc]{$\beta_{3}$}}
\end{picture}
\caption{Illustration of the proof of Corollary~\ref{c3}.}\label{fig:f3}
\end{figure}
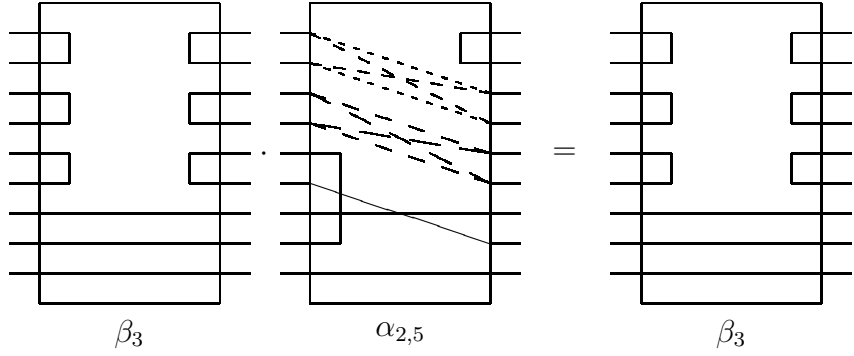

\begin{proof}
Assume that $j$ and $n$ have the same parity, set $k=\frac{n-j+2}{2}$, 
and consider the product $\gamma=\beta_{k}\alpha_{i,j}$. We claim that
for any $s\in\{1,\dots,j-2\}$ there exists $t\in M'$ such that 
$s\equiv_{\gamma}t'$ (see illustration on Figure~\ref{fig:f3}). 
Indeed, for $s\neq i$ this follows directly
from the definitions of $\beta_k$ and $\alpha_{i,j}$. For $s=i$
we have $i\equiv_{\beta_k}i'$, $i\equiv_{\alpha_{i,j}}j$,
$j'\equiv_{\beta_k}(j-1)'$ and $(j-1)\equiv_{\alpha_{i,j}}t'$
for some $t\in M$, which is exactly what we wanted to prove. This
implies that $\gamma\mathcal{R}\beta_{k}$ and hence $\gamma=\beta_{k}$
since $\Lambda$ is an $\mathcal{R}$--cross-section. Now all
equalities in \eqref{c3.1} follow from $\gamma=\beta_{k}$ and the
definitions of $\beta_{k}$ and $\alpha_{i,j}$. One also proves
\eqref{c3.2} by analogous arguments.

To prove \eqref{c3.3} we use induction on $k$. For $k=1$ we
consider $\beta_{1}\alpha_{i,j}$. We obviously have
$\mathrm{corank}(\beta_{1}\alpha_{i,j})=4$ and as $\Lambda$ is
canonical, we get that either $n\equiv_{\alpha_{i,j}}(n-2)'$
and $(n-1)\equiv_{\alpha_{i,j}}(n-3)'$, or else
$n\equiv_{\alpha_{i,j}}(n-3)'$ and $(n-1)\equiv_{\alpha_{i,j}}(n-2)'$.
This implies our statement for $k=1$. Using analogous
arguments we consider $\beta_{2}\alpha_{i,j}$ and proceed by induction.
This completes the proof.
\end{proof}

\begin{proposition}\label{prnew1}
Let $k\in\{1,\dots,\lfloor\frac{n}{2}\rfloor\}$.
\begin{enumerate}[(i)]
\item\label{prnew1.1} Let $\alpha\in \Lambda_k$ and $\eta\in\Lambda_1$.
Then $\alpha\eta=\alpha$ if $\eta\not\in\Lambda_1^{(k+1)}$ and
$\alpha\eta\in \Lambda_{k+1}$ otherwise.
\item\label{prnew1.2} Let $\eta_1,\dots,\eta_k\in\Lambda_1$. Then $\eta_1\cdot\dots\cdot\eta_k\in\Lambda_k$
if and only if $\eta_i\in \Lambda_1^{(i)}$ for all $i\in 
\{1,\dots,k\}$.
\item\label{prnew1.3} If $k<\lfloor\frac{n}{2}\rfloor$ then
$\Lambda_k\Lambda_1^{(k+1)}=\Lambda_{k+1}$. In particular, 
$\Lambda$ is generated by $\Lambda_1$ as a monoid.
\end{enumerate}
\end{proposition}

\begin{proof}
We start with \eqref{prnew1.1}. From the definition of a canonical
$\mathcal{R}$--cross-section we have
\begin{displaymath}
n'\equiv_{\a}(n-1)',(n-2)'\equiv_{\a}(n-3)',\dots, 
(n-2(k-1))'\equiv_{\a}(n-2(k-1)-1)',
\end{displaymath}
and for every $j<n-2(k-1)-1$ there exists $x\in M$ such that
$x\equiv_{\a}j'$. If $\gamma\not\in \Lambda_1^{(k+1)}$ then $\rank(\alpha\gamma)=\rank(\alpha)$. Thus $\alpha\mathcal{R} 
\alpha\gamma$ and hence $\alpha=\alpha\gamma$. In the case when  
$\gamma\in \Lambda_1^{(k+1)}$ one gets 
$\rank(\alpha\gamma)=\rank(\alpha)-2$, which proves \eqref{prnew1.1}.
\eqref{prnew1.2} follows immediately from \eqref{prnew1.1} by
induction on $k$.

Let $k<\lfloor\frac{n}{2}\rfloor$ and $\a\in \Lambda_{k+1}$.
Let further $X_1,\dots,X_{k+1}$ be disjoint two-element subsets
of $M$ contained in $\a$. Let $\beta$ be the unique element of
$\Lambda_k$, containing $X_1,\dots,X_{k}$. Let $X_{k+1}=\{a,b\}$
and  $u,v\in M$ be such that $a\equiv_{\beta}u'$ and
$b\equiv_{\beta}v'$. Let $i<j$ be such that $\{i,j\}=\{u,v\}$.
We have $\beta\alpha_{i,j}\mathcal{R}\a$ by
construction, implying $\beta\alpha_{i,j}=\a$ since
$\Lambda$ is an $\mathcal{R}$--cross-section. This proves
\eqref{prnew1.3}.
\end{proof}

\begin{lemma}\label{l5}
For all $k=1,\dots,\lfloor\frac{n}{2}\rfloor$ we have 
$|\Lambda_k|=\binom{n}{2k}(2k-1)!!$.
\end{lemma}

\begin{proof}
Let $k\in\{1,\dots,\lfloor\frac{n}{2}\rfloor\}$. Then $|\Lambda_k|$
coincides with the number of $\mathcal{R}$--classes of corank $2k$. To define 
such an $\mathcal{R}$--class we have to choose $k$ subsets of $M$, each
containing $2$ elements, when the order of subsets is not important.
This can be done in $\binom{n}{2k}|\B_k|=\binom{n}{2k}(2k-1)!!$ different ways,
completing the proof.
\end{proof}

The following statement is the key observation in our attempt to understand the
structure of $\Lambda$.

\begin{proposition}\label{decomposition}
Let $k\in\{1,\dots,\lfloor\frac{n}{2}\rfloor\}$ and $\a\in\Lambda_k$.
Then the pre-image of $\a$ under the multiplication map
\begin{displaymath}
\begin{array}{rccc}
\mathrm{mult:} & \Lambda_1^{(1)}\times\Lambda_1^{(2)}\times\dots\times\Lambda_1^{(k)}
& \rightarrow & \Lambda_k\\
&(\eta_1,\eta_2,\dots,\eta_k)&\mapsto& \eta_1\eta_2\cdot\dots\cdot\eta_k 
\end{array}
\end{displaymath}
consists of exactly $k!$ elements.
\end{proposition}

\begin{proof}
The map $\mathrm{mult}$ is well-defined by Proposition~\ref{prnew1}\eqref{prnew1.2}
and is surjective by Proposition~\ref{prnew1}\eqref{prnew1.3}. Let
$X_1,\dots,X_k$ be disjoint two-element subsets of $M$, contained in $\a$,
and let $\sigma\in\mathcal{S}_k$. For $i\in\{1,\dots,k\}$ denote by
$\gamma_i^{(\sigma)}$ the unique element of $\Lambda_i$ containing
$X_{\sigma(1)},X_{\sigma(2)},\dots,X_{\sigma(i)}$. Set for convenience
$\gamma_0^{(\sigma)}=\mathrm{id}$.  Then the same arguments as
in the proof of Proposition~\ref{prnew1}\eqref{prnew1.3} give a unique 
element, $\eta_i^{(\sigma)}\in\Lambda_1^{(1)}$, such that 
$\gamma_{i-1}^{(\sigma)}\eta_i^{(\sigma)}=\gamma_i^{(\sigma)}$. By construction
we have $\eta_1^{(\sigma)}\eta_2^{(\sigma)}\cdot\dots\cdot\eta_k^{(\sigma)}
\mathcal{R}\a$ and hence
$\mathrm{mult}(\eta_1^{(\sigma)},\eta_2^{(\sigma)},\dots,\eta_k^{(\sigma)})=\a$.
Moreover, obviously,
\begin{displaymath}
(\eta_1^{(\sigma)},\eta_2^{(\sigma)},\dots,\eta_k^{(\sigma)})\neq
(\eta_1^{(\sigma')},\eta_2^{(\sigma')},\dots,\eta_k^{(\sigma')})
\end{displaymath}
if $\sigma\neq\sigma'$. This means that the pre-image of
$\a$ under $\mathrm{mult}$ consists of at least $k!$ elements.
From the surjectivity of $\mathrm{mult}$ it follows that
\begin{equation}\label{eqcomn}
|\Lambda_1^{(1)}\times\Lambda_1^{(2)}\times\dots\times\Lambda_1^{(k)}|\geq
|\Lambda_k|\cdot k!.
\end{equation}
However, we have 
\begin{equation}\label{eqcomb}
|\Lambda_1^{(1)}\times\Lambda_1^{(2)}\times\dots\times\Lambda_1^{(k)}|=
\binom{n}{2}\binom{n-2}{2}\dots\binom{n-2(k-1)}{2}
\end{equation}
and, comparing \eqref{eqcomn}, \eqref{eqcomb} and Lemma~\ref{l5}, we
obtain that \eqref{eqcomn} must be an equality. The statement follows.
\end{proof}

We also have the following recursion for the canonical 
$\mathcal{R}$--cross-section. 

\begin{proposition}\label{prnewrec}
Let $n>2$ and $\Lambda$ be a canonical $\mathcal{R}$--cross-section of $\B_n$.
\begin{enumerate}[(i)]
\item\label{prnewrec.1}
For every $\beta\in \Lambda$ the element 
$(\alpha_{n-1,n}\beta)|_{\{1,\dots,n-2\}}$
is a well-defined element of $\B_{n-2}$ and the set
$\Phi=\{(\alpha_{n-1,n}\beta)|_{\{1,\dots,n-2\}}\,:\,\beta\in \Lambda\}$ is
a canonical $\mathcal{R}$--cross-section of $\B_{n-2}$.
\item\label{prnewrec.2} The map
\begin{displaymath}
\begin{array}{rccc}
\varphi: & \Lambda & \rightarrow & \Phi\\
& \beta & \mapsto & \alpha_{n-1,n}\beta
\end{array}
\end{displaymath}
is a homomorphism, which sends idempotent elements of corank $2$ to the
identity and is injective on the set of all nonidempotent elements of corank $2$.
\end{enumerate}
\end{proposition}

\begin{proof}
Consider the map $\overline{\varphi}:\Lambda\to\Lambda$, defined
via $\overline{\varphi}(\beta)=\alpha_{n-1,n}\beta\alpha_{n-1,n}$
for all $\beta\in \Lambda$. A direct calculation shows that 
$\overline{\varphi}(\beta)=\alpha_{n-1,n}\beta$, which implies that
$\overline{\varphi}$ is a homomorphism. It is easy to see that the
image of $\overline{\varphi}$ coincides with
$N=\{\beta\in\Lambda\,:\,\{n-1,n\},\{(n-1)',n'\}\in\beta\}$.
In particular, $\{1,\dots,n-2\}$ is $\beta$-invariant for every
$\beta\in N$ and hence for all such $\beta$ the element
$\beta|_{\{1,\dots,n-2\}}$ is a well-defined element of $\B_{n-2}$.
Hence $\varphi$ is well-defined. Moreover, $\varphi$ is a homomorphism
since both $\overline{\varphi}:\Lambda\to N$ and
${}_-|_{\{1,\dots,n-2\}}:N\to\B_{n-2}$ are. This proves
the first parts of \eqref{prnewrec.1} and \eqref{prnewrec.2}.
The rest of \eqref{prnewrec.2} is proved by a direct calculation.

So, we are left to show that $\varphi(\Lambda)$ is a
canonical $\mathcal{R}$--cross-section of $\B_{n-2}$.
Forgetting $\{n-1,n\}$ identifies the elements in the set
of all collections of two-element subsets of
$\{1,\dots,n\}$, containing $\{n-1,n\}$, and the set
 of all collections of two-element subsets of
$\{1,\dots,n-2\}$. Since $\Lambda$ was an 
$\mathcal{R}$--cross-section of $\B_{n}$, it follows that
the map ${}_-|_{\{1,\dots,n-2\}}$ gives rise to a bijection
between the elements of $N$ and all collections of two-element 
subsets of $\{1,\dots,n-2\}$. This means that 
$\varphi(\Lambda)=N|_{\{1,\dots,n-2\}}$ is an 
$\mathcal{R}$--cross-section of $\B_{n-2}$. It is easy to see
that this cross-section is canonical. This completes the proof.
\end{proof}

Proposition~\ref{decomposition} implies that $\Lambda$ is completely determined 
(generated) by $\Lambda_1$. Recall that we are still working with a given
$\Lambda$ (and we still do not know if it exists). $\Lambda_1$ consists of
$\a_{i,j}$ and hence to describe $\Lambda_1$ we have to determine all
$\a_{i,j}$ explicitly. If $j\in\{n-1,n\}$ then $\a_{i,j}$ is an idempotent by
Lemma~\ref{l1}\eqref{l1.2} and it is explicitly described by
Lemma~\ref{l1}\eqref{l1.n3}. In all other cases Corollary~\ref{c3}
gives only a precise description of some part of $\a_{i,j}$, since the statement
of Corollary~\ref{c3}\eqref{c3.3} describes certain parts of $\a_{i,j}$, namely the
partition sets, containing $a\in M$ for $a>j$ or $a>j+1$ depending on some parities, 
only up to a bijection between two two-element sets (that is, roughly speaking, up 
to an element of $\S_2$). Our idea now is to write these undetermined parts as 
``parameters'', identifying each of the parameters with an element of $\S_2$
(since in some sense they behave well under multiplication, see Lemma~\ref{l8} below), 
and to investigate the  relations between these parameters. Let $1\leq i<j<n-1$. For
$l=1,\dots,\lfloor\frac{n-j}{2}\rfloor$ we define $\alpha_{i,j}^{(l)}\in \S_2$  in
the following way:
\begin{equation}\label{eq2}
\alpha_{i,j}^{(l)}=
\begin{cases}
\mathrm{id}, & (n-2l+2)\equiv_{\alpha_{i,j}}(n-2l)'\\
(1,2), & (n-2l+2)\equiv_{\alpha_{i,j}}(n-2l-1)'
\end{cases}.
\end{equation}
This definition is motivated by the following easy observation.

\begin{lemma}\label{l8}
Let $1\leq i<j<n-1$, $l\in \{1,\dots,\lfloor\frac{n-j}{2}\rfloor\}$,
and $1\leq s<t<n-1$ be such that $\lfloor\frac{n-t}{2}\rfloor\geq l+1$. 
Let $\beta=\alpha_{i,j}\alpha_{s,t}$. 
\begin{enumerate}[(i)]
\item\label{l8.1} For every $u\in\{n-2l+1,n-2l+2\}$ there exists
$v\in \{n-2l-3,n-2l-2\}$ such that $u\equiv_{\beta}v$.
\item\label{l8.2} Define $\beta^{(l)}\in \S_2$ as follows:
\begin{equation}\label{eq3}
\beta^{(l)}=
\begin{cases}
\mathrm{id}, & (n-2l+1)\equiv_{\beta}(n-2l-3)'\\
(1,2), & (n-2l+1)\equiv_{\beta}(n-2l-2)'
\end{cases}.
\end{equation}
Then $\beta^{(l)}=\alpha_{i,j}^{(l)}\alpha_{s,t}^{(l+1)}$.
\end{enumerate}
\end{lemma}

\begin{proof}
\eqref{l8.1} follows from Corollary~\ref{c3}\eqref{c3.3}, and
\eqref{l8.2} follows from the definitions \eqref{eq2} and \eqref{eq3}.
\end{proof}

Now we would like to describe the canonical $\mathcal{R}$--cross-sections for
small values of $n$.

\begin{proposition}\label{p10}
\begin{enumerate}[(i)]
\item\label{p10.1} For $n=1$ we have $1$ trivial $\mathcal{R}$--cross-section.
\item\label{p10.2} For $n=2$ we have only $1$ $\mathcal{R}$--cross-section, moreover,
it is canonical and consists of $\alpha_{1,2}$ and $\mathrm{id}$.
\item\label{p10.3} For $n=3$ we have $1$ canonical $\mathcal{R}$--cross-section,
consisting of idempotents $\alpha_{1,2}$, $\alpha_{2,3}$, $\alpha_{1,3}$ and $\mathrm{id}$.
\item\label{p10.4} For $n=4$ we have $2$ canonical $\mathcal{R}$--cross-section,
for one of them we have $\alpha_{1,2}^{(1)}=\mathrm{id}$, for another one we
have $\alpha_{1,2}^{(1)}=(1,2)$.
\item\label{p10.5} For $n=5$ we have $8$ canonical $\mathcal{R}$--cross-section,
which correspond to independent choices of the parameters
$\alpha_{1,2}^{(1)},\alpha_{2,3}^{(1)},\alpha_{1,3}^{(1)}\in \S_2$.
\item\label{p10.6} For $n=6$ we have $16$ canonical $\mathcal{R}$--cross-section, 
which correspond to independent choices of the parameters
$\alpha_{1,2}^{(1)},\alpha_{1,2}^{(2)},\alpha_{2,3}^{(1)},\alpha_{1,3}^{(1)}\in \S_2$.
\end{enumerate}
\end{proposition}

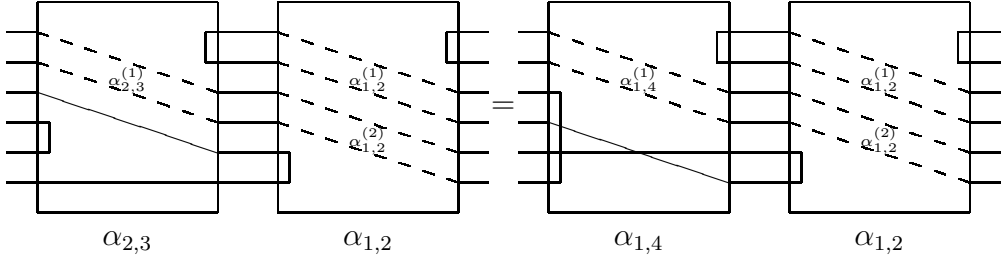
\begin{figure}
\special{em:linewidth 0.4pt}
\unitlength 0.80mm
\linethickness{0.4pt}
\begin{picture}(200.00,50.00)
\drawline(05.00,10.00)(05.00,45.00)
\drawline(05.00,10.00)(35.00,10.00)
\drawline(05.00,45.00)(35.00,45.00)
\drawline(35.00,10.00)(35.00,45.00)
\drawline(00.00,15.00)(05.00,15.00)
\drawline(00.00,20.00)(05.00,20.00)
\drawline(00.00,25.00)(05.00,25.00)
\drawline(00.00,30.00)(05.00,30.00)
\drawline(00.00,35.00)(05.00,35.00)
\drawline(00.00,40.00)(05.00,40.00)
\drawline(40.00,15.00)(35.00,15.00)
\drawline(40.00,20.00)(35.00,20.00)
\drawline(40.00,25.00)(35.00,25.00)
\drawline(40.00,30.00)(35.00,30.00)
\drawline(40.00,35.00)(35.00,35.00)
\drawline(40.00,40.00)(35.00,40.00)
\drawline(05.00,15.00)(35.00,15.00)
\drawline(05.00,30.00)(35.00,20.00)
\drawline(05.00,20.00)(07.00,20.00)
\drawline(05.00,25.00)(07.00,25.00)
\drawline(07.00,20.00)(07.00,25.00)
\drawline(35.00,40.00)(33.00,40.00)
\drawline(35.00,35.00)(33.00,35.00)
\drawline(33.00,40.00)(33.00,35.00)
\dashline{2}(05.00,40.00)(35.00,30.00)
\dashline{2}(05.00,35.00)(35.00,25.00)
\drawline(45.00,10.00)(45.00,45.00)
\drawline(45.00,10.00)(75.00,10.00)
\drawline(45.00,45.00)(75.00,45.00)
\drawline(75.00,10.00)(75.00,45.00)
\drawline(40.00,15.00)(45.00,15.00)
\drawline(40.00,20.00)(45.00,20.00)
\drawline(40.00,25.00)(45.00,25.00)
\drawline(40.00,30.00)(45.00,30.00)
\drawline(40.00,35.00)(45.00,35.00)
\drawline(40.00,40.00)(45.00,40.00)
\drawline(80.00,15.00)(75.00,15.00)
\drawline(80.00,20.00)(75.00,20.00)
\drawline(80.00,25.00)(75.00,25.00)
\drawline(80.00,30.00)(75.00,30.00)
\drawline(80.00,35.00)(75.00,35.00)
\drawline(80.00,40.00)(75.00,40.00)
\dashline{2}(45.00,25.00)(75.00,15.00)
\dashline{2}(45.00,30.00)(75.00,20.00)
\drawline(45.00,15.00)(47.00,15.00)
\drawline(45.00,20.00)(47.00,20.00)
\drawline(47.00,15.00)(47.00,20.00)
\drawline(75.00,40.00)(73.00,40.00)
\drawline(75.00,35.00)(73.00,35.00)
\drawline(73.00,40.00)(73.00,35.00)
\dashline{2}(45.00,40.00)(75.00,30.00)
\dashline{2}(45.00,35.00)(75.00,25.00)
\drawline(90.00,10.00)(90.00,45.00)
\drawline(90.00,10.00)(120.00,10.00)
\drawline(90.00,45.00)(120.00,45.00)
\drawline(120.00,10.00)(120.00,45.00)
\drawline(85.00,15.00)(90.00,15.00)
\drawline(85.00,20.00)(90.00,20.00)
\drawline(85.00,25.00)(90.00,25.00)
\drawline(85.00,30.00)(90.00,30.00)
\drawline(85.00,35.00)(90.00,35.00)
\drawline(85.00,40.00)(90.00,40.00)
\drawline(125.00,15.00)(120.00,15.00)
\drawline(125.00,20.00)(120.00,20.00)
\drawline(125.00,25.00)(120.00,25.00)
\drawline(125.00,30.00)(120.00,30.00)
\drawline(125.00,35.00)(120.00,35.00)
\drawline(125.00,40.00)(120.00,40.00)
\drawline(90.00,20.00)(120.00,20.00)
\drawline(90.00,25.00)(120.00,15.00)
\drawline(90.00,15.00)(92.00,15.00)
\drawline(90.00,30.00)(92.00,30.00)
\drawline(92.00,15.00)(92.00,30.00)
\drawline(120.00,40.00)(118.00,40.00)
\drawline(120.00,35.00)(118.00,35.00)
\drawline(118.00,40.00)(118.00,35.00)
\dashline{2}(90.00,40.00)(120.00,30.00)
\dashline{2}(90.00,35.00)(120.00,25.00)
\drawline(130.00,10.00)(130.00,45.00)
\drawline(130.00,10.00)(160.00,10.00)
\drawline(130.00,45.00)(160.00,45.00)
\drawline(160.00,10.00)(160.00,45.00)
\drawline(125.00,15.00)(130.00,15.00)
\drawline(125.00,20.00)(130.00,20.00)
\drawline(125.00,25.00)(130.00,25.00)
\drawline(125.00,30.00)(130.00,30.00)
\drawline(125.00,35.00)(130.00,35.00)
\drawline(125.00,40.00)(130.00,40.00)
\drawline(165.00,15.00)(160.00,15.00)
\drawline(165.00,20.00)(160.00,20.00)
\drawline(165.00,25.00)(160.00,25.00)
\drawline(165.00,30.00)(160.00,30.00)
\drawline(165.00,35.00)(160.00,35.00)
\drawline(165.00,40.00)(160.00,40.00)
\dashline{2}(130.00,25.00)(160.00,15.00)
\dashline{2}(130.00,30.00)(160.00,20.00)
\drawline(130.00,15.00)(132.00,15.00)
\drawline(130.00,20.00)(132.00,20.00)
\drawline(132.00,15.00)(132.00,20.00)
\drawline(160.00,40.00)(158.00,40.00)
\drawline(160.00,35.00)(158.00,35.00)
\drawline(158.00,40.00)(158.00,35.00)
\dashline{2}(130.00,40.00)(160.00,30.00)
\dashline{2}(130.00,35.00)(160.00,25.00)
\put(82.50,27.50){\makebox(0,0)[cc]{$=$}}
\put(20.00,5.00){\makebox(0,0)[cc]{$\alpha_{2,3}$}}
\put(20.00,32.00){\makebox(0,0)[cc]{\tiny $\alpha_{2,3}^{(1)}$}}
\put(60.00,32.00){\makebox(0,0)[cc]{\tiny $\alpha_{1,2}^{(1)}$}}
\put(60.00,22.00){\makebox(0,0)[cc]{\tiny $\alpha_{1,2}^{(2)}$}}
\put(105.00,32.00){\makebox(0,0)[cc]{\tiny $\alpha_{1,4}^{(1)}$}}
\put(145.00,32.00){\makebox(0,0)[cc]{\tiny $\alpha_{1,2}^{(1)}$}}
\put(145.00,22.00){\makebox(0,0)[cc]{\tiny $\alpha_{1,2}^{(2)}$}}
\put(60.00,5.00){\makebox(0,0)[cc]{$\alpha_{1,2}$}}
\put(105.00,5.00){\makebox(0,0)[cc]{$\alpha_{1,4}$}}
\put(145.00,5.00){\makebox(0,0)[cc]{$\alpha_{1,2}$}}
\end{picture}
\caption{Illustration of the equality \eqref{relp10.2} implying $\alpha_{2,3}^{(1)}=\alpha_{1,4}^{(1)}$.}\label{fig:f4}
\end{figure}

\begin{proof}
The statements \eqref{p10.1}, \eqref{p10.2}, and \eqref{p10.3} are obvious. For
$n=4$ we observe that there is only one parameter, namely $\alpha_{1,2}^{(1)}\in\S_2$.
A direct calculation shows that both values of the parameter indeed lead to
cross-sections. This proves \eqref{p10.4}. For $n=5$ we observe that there are exactly
$3$ parameters, namely $\alpha_{1,2}^{(1)}\in\S_2$, $\alpha_{2,3}^{(1)}\in\S_2$, and 
$\alpha_{1,3}^{(1)}\in\S_2$. A direct calculation again shows that all values of these
parameters indeed lead to cross-sections. This proves \eqref{p10.5}.

Let us now consider the case $n=6$. In this case we have $7$ parameters, namely
$\alpha_{1,2}^{(1)}\in\S_2$, $\alpha_{1,2}^{(2)}\in\S_2$, $\alpha_{1,3}^{(1)}\in\S_2$, $\alpha_{1,4}^{(1)}\in\S_2$, $\alpha_{2,3}^{(1)}\in\S_2$, $\alpha_{2,4}^{(1)}\in\S_2$, and 
$\alpha_{3,4}^{(1)}\in\S_2$. Using the fact that $\Lambda$ is an 
$\mathcal{R}$--cross-section, by a direct calculation we obtain that the following 
relations should be satisfied (see Figure~\ref{fig:f4}):
\begin{gather}
\alpha_{1,3}\alpha_{1,2}=\alpha_{2,4}\alpha_{1,2},\label{relp10.1}\\
\alpha_{2,3}\alpha_{1,2}=\alpha_{1,4}\alpha_{1,2},\label{relp10.2}\\
\alpha_{1,2}\alpha_{1,2}=\alpha_{3,4}\alpha_{1,2}.\label{relp10.3}
\end{gather}
Using now Lemma~\ref{l8}, from \eqref{relp10.1} we obtain 
$\alpha_{1,3}^{(1)}\alpha_{1,2}^{(2)}=\alpha_{2,4}^{(1)}\alpha_{1,2}^{(2)}$
implying $\alpha_{1,3}^{(1)}=\alpha_{2,4}^{(1)}$. Analogously 
from \eqref{relp10.2} we obtain $\alpha_{2,3}^{(1)}=\alpha_{1,4}^{(1)}$,
from \eqref{relp10.3} we obtain $\alpha_{1,2}^{(1)}=\alpha_{3,4}^{(1)}$. This implies that
all parameters can be expressed in terms of $\alpha_{1,2}^{(1)}$, $\alpha_{1,2}^{(2)}$,
$\alpha_{2,3}^{(1)}$, and $\alpha_{1,3}^{(1)}$. A direct (but quite long) calculation 
shows that all values of these parameters indeed lead to distinct
cross-sections. This proves \eqref{p10.6}.
\end{proof}

Later on we will also need the following relation between the 
(canonical) $\mathcal{R}$--cross-sec\-ti\-ons for different $n$.

\begin{proposition}\label{p4}
Let $\Gamma=\{\alpha\in\Lambda\,:\, 1\equiv_{\a}1'\}$. Then the set
$\overline{\Gamma}=\{\alpha|_{\{2,\dots,n\}}\,:\a\in \Gamma\}$  is an
$\mathcal{R}$--cross-section of $\B(\{2,\dots,n\})$. Identifying
$\{2,\dots,n\}$ with $\{1,\dots,n-1\}$ via $x\mapsto x-1$, 
$\overline{\Gamma}$ becomes a canonical $\mathcal{R}$--cross-section of 
$\B_{n-1}$.
\end{proposition}

\begin{proof}
Left to the reader.
\end{proof}

Now we are ready to go to the general case. Denote by $\mathfrak{A}_n$ the
set of all pairs $(n-3,n-2)$, $(n-5,n-4)$,\dots.

\begin{proposition}\label{p11}
\begin{enumerate}[(i)]
\item \label{p5.1} Let $n\geq 8$ be even, then we have the following equalities:
\begin{gather}
\alpha_{1,2}^{(l)}=\alpha_{i,j}^{(l)}\quad\text{ for all }\quad
(i,j)\in\mathfrak{A}_n \text{ and for all } l=1,\dots, \frac{n-j}{2}\label{relp11.1},\\
\alpha_{2,3}^{(l)}=\alpha_{i,j}^{(l)}\quad\text{ for all }\quad
(i,j)\not\in\mathfrak{A}_n \text{ and for all } l=1,\dots, \lfloor\frac{n-j}{2}\rfloor
\label{relp11.2},\\
\alpha_{1,2}^{(l)}\alpha_{2,3}^{(l+1)}=\alpha_{2,3}^{(l)}\alpha_{1,2}^{(l+1)}
\text{ for all }l=1,\dots, \frac{n}{2}-3\label{relp11.25}.
\end{gather}
\item \label{p5.2} Let $n\geq 7$ be odd, then we have the following equalities:
\begin{gather}
\alpha_{2,3}^{(l)}=\alpha_{i,j}^{(l)}\quad\text{ for all }\quad
(i,j)\in\mathfrak{A}_n \text{ and for all } l=1,\dots, \frac{n-j}{2}\label{relp11.3},\\
\alpha_{1,2}^{(l)}=\alpha_{i,j}^{(l)}\quad\text{ for all }\quad
(i,j)\not\in\mathfrak{A}_n \text{ and for all } l=1,\dots, \lfloor\frac{n-j}{2}\rfloor
\label{relp11.4},\\
\alpha_{1,2}^{(l)}\alpha_{2,3}^{(l+1)}=\alpha_{2,3}^{(l)}\alpha_{1,2}^{(l+1)}
\text{ for all }l=1,\dots, \frac{n+1}{2}-3\label{relp11.45}.
\end{gather}
\end{enumerate}
\end{proposition}

\begin{proof}
We prove both statements using induction on $n$. We start with $n=7$. In this case we use the
fact that $\Lambda$ is an $\mathcal{R}$--cross-section to obtain the following equalities:
\begin{gather}
\label{relp11.p1} \alpha_{2,3}\alpha_{2,3}=\alpha_{4,5}\alpha_{2,3},\\
\label{relp11.p2} \alpha_{2,4}\alpha_{1,2}=\alpha_{1,5}\alpha_{1,2},\\
\label{relp11.p3} \alpha_{1,4}\alpha_{1,2}=\alpha_{2,5}\alpha_{1,2},\\
\label{relp11.p4} \alpha_{3,5}\alpha_{2,3}=\alpha_{2,4}\alpha_{2,3},\\
\label{relp11.p5} \alpha_{3,4}\alpha_{2,3}=\alpha_{2,5}\alpha_{2,3},\\
\label{relp11.p6} \alpha_{1,5}\alpha_{1,3}=\alpha_{3,4}\alpha_{1,3},\\
\label{relp11.p7} \alpha_{1,3}\alpha_{2,3}=\alpha_{4,5}\alpha_{1,3}.
\end{gather}
The arguments, analogous to those in the proof of Proposition~\ref{p10} give the following:
\eqref{relp11.p1} implies $\alpha_{2,3}^{(1)}=\alpha_{4,5}^{(1)}$. Further, 
\eqref{relp11.p2} implies $\alpha_{2,4}^{(1)}=\alpha_{1,5}^{(1)}$, 
\eqref{relp11.p3} implies $\alpha_{1,4}^{(1)}=\alpha_{2,5}^{(1)}$, 
\eqref{relp11.p4} implies $\alpha_{3,5}^{(1)}=\alpha_{2,4}^{(1)}$, 
\eqref{relp11.p5} implies $\alpha_{3,4}^{(1)}=\alpha_{2,5}^{(1)}$,
\eqref{relp11.p6} implies $\alpha_{1,5}^{(1)}=\alpha_{3,4}^{(1)}$. This implies
\begin{displaymath}
\alpha_{1,5}^{(1)}=\alpha_{1,4}^{(1)}=\alpha_{2,4}^{(1)}=\alpha_{2,5}^{(1)}=
\alpha_{3,4}^{(1)}=\alpha_{3,5}^{(1)}.
\end{displaymath}
Moreover,  \eqref{relp11.p7} implies 
\begin{equation}\label{relp11.p8}
\alpha_{1,3}^{(1)}\alpha_{2,3}^{(2)}=\alpha_{4,5}^{(1)}\alpha_{1,3}^{(2)}.
\end{equation}
Now we have to go to the case-by-case analysis. We consider two cases,
$\alpha_{1,3}^{(2)}=\mathrm{id}$ and $\alpha_{1,3}^{(2)}=(1,2)$. 

For $\alpha_{1,3}^{(2)}=\mathrm{id}$ we have the identities
$\alpha_{1,3}\alpha_{1,2}=\alpha_{2,4}\alpha_{1,3}$ and
$\alpha_{1,3}\alpha_{1,3}=\alpha_{2,5}\alpha_{1,3}$ giving
\begin{equation}\label{relp11.p9}
\alpha_{1,3}^{(1)}\alpha_{1,2}^{(2)}=\alpha_{2,4}^{(1)}\alpha_{1,3}^{(2)}
\text{ and }
\alpha_{1,3}^{(1)}\alpha_{1,3}^{(2)}=\alpha_{2,5}^{(1)}\alpha_{1,3}^{(2)},
\end{equation}
respectively.

For $\alpha_{1,3}^{(2)}=(1,2)$ we have the identities
$\alpha_{1,3}\alpha_{1,2}=\alpha_{2,5}\alpha_{1,3}$ and
$\alpha_{1,3}\alpha_{1,3}=\alpha_{2,4}\alpha_{1,3}$ giving
\begin{equation}\label{relp11.p10}
\alpha_{1,3}^{(1)}\alpha_{1,2}^{(2)}=\alpha_{2,5}^{(1)}\alpha_{1,3}^{(2)}
\text{ and }
\alpha_{1,3}^{(1)}\alpha_{1,3}^{(2)}=\alpha_{2,4}^{(1)}\alpha_{1,3}^{(2)},
\end{equation}
respectively. Since we already know that $\alpha_{2,5}^{(1)}=\alpha_{2,4}^{(1)}$,
we have that both \eqref{relp11.p9} and \eqref{relp11.p10} in  fact do not
depend on the values of $\alpha_{1,3}^{(2)}$.

Combining \eqref{relp11.p8} and \eqref{relp11.p10} we consequently obtain
$\alpha_{1,3}^{(1)}=\alpha_{2,5}^{(1)}$, $\alpha_{1,3}^{(2)}=\alpha_{1,2}^{(2)}$
and $\alpha_{1,3}^{(1)}\alpha_{2,3}^{(2)}=\alpha_{4,5}^{(1)}\alpha_{1,2}^{(2)}$. 
To complete the proof we are now left to show that 
$\alpha_{1,3}^{(1)}=\alpha_{1,2}^{(1)}$

Consider two cases: $\alpha_{1,2}^{(2)}=\mathrm{id}$ and 
$\alpha_{1,2}^{(2)}=(1,2)$. In the first case we obtain 
$\alpha_{1,2}\alpha_{1,2}=\alpha_{3,4}\alpha_{1,2}$ implying
\begin{equation}\label{relp11.p11}
\alpha_{1,2}^{(1)}\alpha_{1,2}^{(2)}=\alpha_{3,4}^{(1)}\alpha_{1,2}^{(2)}.
\end{equation}
In the second case case we obtain 
$\alpha_{1,2}\alpha_{1,2}=\alpha_{3,5}\alpha_{1,2}$ implying
$\alpha_{1,2}^{(1)}\alpha_{1,2}^{(2)}=\alpha_{3,5}^{(1)}\alpha_{1,2}^{(2)}$,
which is the same as \eqref{relp11.p11} as we already know that
$\alpha_{3,4}^{(1)}=\alpha_{3,5}^{(1)}$. This implies that
\eqref{relp11.p11} holds in all cases, which gives
$\alpha_{1,2}^{(1)}=\alpha_{3,4}^{(1)}=\alpha_{2,5}^{(1)}=\alpha_{1,3}^{(1)}$. 
So, the case $n=7$ is complete.

Now we prove the induction step and, because of the inductive assumption
and Proposition~\ref{p4}, it is enough to prove either the statement \eqref{relp11.4}
or the statements \eqref{relp11.1}, \eqref{relp11.2}, respectively, for the elements 
$\alpha_{1,j}$, $j=2,\dots,n-2$, and either the statement \eqref{relp11.25} or 
\eqref{relp11.45}, respectively, depending on the parity of $n$. 

Assume that $n$ is odd. Then $\alpha_{1,j}\not\in \mathfrak{A}_n$ for all $j=2,\dots,n-2$
and hence we have to check \eqref{relp11.4}. For even $j\geq 4$ we have 
$\alpha_{1,j}\alpha_{1,2}=\alpha_{2,j+1}\alpha_{1,2}$ and for odd $j\geq 5$ we have
$\alpha_{1,j}\alpha_{1,2}=\alpha_{2,j-1}\alpha_{1,2}$ and in both cases we obtain
$\alpha_{1,j}^{(l)}=\alpha_{3,4}^{(l)}$ for all possible $l$ by inductive assumptions.

Let $l=\frac{n-3}{2}$. Then for $\alpha_{1,2}^{(l)}=\mathrm{id}$ we have
$\alpha_{1,2}\alpha_{1,2}=\alpha_{3,4}\alpha_{1,2}$ and for 
$\alpha_{1,2}^{(l)}=(1,2)$ we have
$\alpha_{1,2}\alpha_{1,2}=\alpha_{3,5}\alpha_{1,2}$. But since 
$\alpha_{3,4}^{(s)}=\alpha_{3,5}^{(s)}$ for all possible $s$ by induction,
we obtain $\alpha_{1,2}^{(s)}=\alpha_{3,4}^{(s)}$ for all $s<l$. 

Let $l=\frac{n-3}{2}$. Then for $\alpha_{1,3}^{(l)}=\mathrm{id}$ we have
$\alpha_{1,3}\alpha_{1,3}=\alpha_{2,4}\alpha_{1,3}$ and for 
$\alpha_{1,3}^{(l)}=(1,2)$ we have $\alpha_{1,3}\alpha_{1,3}=\alpha_{2,5}\alpha_{1,3}$. 
But, since $\alpha_{2,4}^{(s)}=\alpha_{2,5}^{(s)}$ for all possible $s$ by induction,
we obtain $\alpha_{1,3}^{(s)}=\alpha_{2,4}^{(s)}$ for all $s<l$. That
$\alpha_{1,3}^{(l)}=\alpha_{1,2}^{(l)}$ is proved using the same arguments as in the
paragraph containing the formula \eqref{relp11.p11}. So, the proof
of \eqref{relp11.4} for the elements $\alpha_{1,j}$, $j=2,\dots,n-2$ is complete.

Finally, we have $\alpha_{1,2}\alpha_{2,3}=\alpha_{4,5}\alpha_{1,2}$ which implies
\eqref{relp11.45}.

Assume now that $n$ is even. Then we have $\alpha_{1,2}\alpha_{1,2}=\alpha_{3,4}\alpha_{1,2}$
which implies \eqref{relp11.1} for $\alpha_{1,2}$. For the elements $\alpha_{1,j}$,
$j\geq 4$ the arguments are the same as in the case of odd $n$. Further, we also have
$\alpha_{1,3}\alpha_{1,2}=\alpha_{2,4}\alpha_{1,2}$ implying equalities
\eqref{relp11.2} for $\alpha_{1,3}$. And finally, \eqref{relp11.25}
follows from the inductive assumption and $\alpha_{2,3}\alpha_{1,2}=\alpha_{1,4}\alpha_{1,2}$.
This completes the proof.
\end{proof}

Now we are ready to construct a canonical 
$\mathcal{R}$--cross-section for $\B_n$
in the general case.

\begin{proposition}\label{p12}
Let $n\in\N$, $l=\lfloor\frac{n-2}{2}\rfloor$, $m=\lfloor\frac{n-3}{2}\rfloor$,
and choose  $x_1,\dots,x_l$, $y_1,\dots,y_m\in\S_2$ such that
$x_iy_{i+1}=y_ix_{i+1}$ for all possible $i$. For $1\leq i<j\leq n$
let $\alpha_{i,j}$ be the element which satisfies the corresponding equalities of
Lemma~\ref{l1}\eqref{l1.n3} and Corollary~\ref{c3}, and, additionally, the following
conditions:
\begin{enumerate}[(a)]
\item $\alpha_{i,j}^{(s)}=x_s$ for all $(i,j)\in \mathfrak{A}_n$ and for all possible $s$;
\item $\alpha_{i,j}^{(s)}=y_s$ for all $(i,j)\not\in \mathfrak{A}_n$ and for all possible $s$.
\end{enumerate}
Then the elements $\mathrm{id}$ and $\alpha_{i,j}$, $1\leq i<j\leq n$, generate a 
canonical $\mathcal{R}$--cross-section of $\B_n$.
\end{proposition}

\begin{proof}
Set $\Gamma=\{\alpha_{i,j}\,:\,1\leq i<j\leq n\}$ and let $\Phi=\langle\Gamma\rangle$
be the monoid, generated by $\Gamma$. For $k=1,2,\dots,\lfloor \frac{n}{2}\rfloor$ set
$\Gamma_k=\{\alpha_{i,j}\,:\,j\leq n-2(k-1)\}$.
Define $\Phi_0=\{\mathrm{id}\}$ and for $i=1,2,\dots,\lfloor \frac{n}{2}\rfloor$ set
$\Phi_i=\Gamma_1\Gamma_2\cdot \dots\cdot\Gamma_i$.
Finally, set $\overline{\Phi}=\cup_{i=0}^{\lfloor \frac{n}{2}\rfloor}\Phi_i$.
We are going to show  that $\overline{\Phi}$ contains exactly
one element of each $\mathcal{R}$--class of $\B_n$, and then that 
$\overline{\Phi}=\Phi$.

By Theorem~\ref{tgreen}\eqref{tgreen.2}, every $\mathcal{R}$--class of $\B_n$ 
of corank $2k$ is uniquely determined by an unordered collection of $k$
disjoint two-element subsets of $M$. Using the same arguments as in 
Proposition~\ref{decomposition} one shows that every $\mathcal{R}$--class of 
$\B_n$ of corank $2k$ contains exactly $k!$ products of the form 
$\eta_1\cdot\dots\cdot\eta_k$, where all $\eta_i\in \Gamma_i$. 

{\em Step 1.} Let us first show that, under the assumptions of our statement,
$\Phi_2$ contains exactly one element of each $\mathcal{R}$--class of corank $4$. 
The arguments above show that every such class contains at most two elements of $\Phi_2$.
Consider $\alpha_{i,j}$ and $\alpha_{s,t}$ such that $t<n-1$. 
Then there exist the unique pair $u<v\in M$ such that
$u\equiv_{\alpha_{i,j}\alpha_{s,t}}v$ and $(u,v)\neq (i,j)$, moreover, there exists
a unique pair, $p<q\in M$, such that 
\begin{equation}\label{brbr1}
\alpha_{i,j}\alpha_{s,t}\mathcal{R}\alpha_{u,v}\alpha_{p,q}.
\end{equation}
It is enough to show that 
\begin{equation}\label{brbr2}
\alpha_{i,j}\alpha_{s,t}=\alpha_{u,v}\alpha_{p,q}. 
\end{equation}

\begin{lemma}\label{lbr1}
\eqref{brbr2} is equivalent to the collection of the following conditions:
\begin{gather}
\alpha_{n-1,n}\alpha_{i,j}\alpha_{s,t}=\alpha_{n-1,n}\alpha_{u,v}\alpha_{p,q}
\label{lbr1.1},\\
(n-1)\equiv_{\alpha_{i,j}\alpha_{s,t}}f\Rightarrow
(n-1)\equiv_{\alpha_{u,v}\alpha_{p,q}}f\quad\text{ for all }\quad f\in M\cup M',
\label{lbr1.2}\\
n\equiv_{\alpha_{i,j}\alpha_{s,t}}f\Rightarrow
n\equiv_{\alpha_{u,v}\alpha_{p,q}}f\quad\text{ for all }\quad f\in M\cup M'.
\label{lbr1.3}
\end{gather}
\end{lemma}

\begin{proof}
That \eqref{brbr2} implies \eqref{lbr1.1}, \eqref{lbr1.2}, and \eqref{lbr1.3}
is obvious. Hence we assume that \eqref{lbr1.1}, \eqref{lbr1.2}, and \eqref{lbr1.3}
are satisfied and we have to prove \eqref{brbr2}. Set
$\a=\alpha_{i,j}\alpha_{s,t}$ and $\beta=\alpha_{u,v}\alpha_{p,q}$. We have
to show that for every $x,y\in M\cup M'$ the condition
$x\equiv_{\a}y$ implies the condition $x\equiv_{\beta}y$. We know that
$\mathrm{corank}(\a)=\mathrm{corank}(\beta)=4$. By the definition
of $\Gamma$ we have that both, $\{(n-1)',n'\}$ and $\{(n-2)',(n-3)'\}$, belong
to both $\a$ and $\beta$. Hence, without loss of generality we can assume
that $x\in M$. If $\{x,y\}\cap\{n-1,n\}\neq \varnothing$, then the necessary
statement follows from \eqref{lbr1.2} and\eqref{lbr1.3}. If
$\{x,y\}\cap\{n-1,n\}=\varnothing$ then it follows from  \eqref{lbr1.1}
since $a\equiv_{\a_{n-1,n}}a'$ for all $a<n-1$.
\end{proof}

Now we claim that \eqref{lbr1.1} follows by induction on $n$ with the cases 
$n=1,2$ being trivial. The idea of the induction is based on the statement 
of Proposition~\ref{prnewrec}.
From the definition of $\alpha_{x,y}$ we have $\alpha_{x,y}\alpha_{n-1,n}=\alpha_{x,y}$
for all appropriate $x,y$. Hence \eqref{lbr1.1} is equivalent to
\begin{equation}\label{brbr2n2}
\alpha_{n-1,n}\alpha_{i,j}\alpha_{n-1,n}\alpha_{s,t}=
\alpha_{n-1,n}\alpha_{u,v}\alpha_{n-1,n}\alpha_{p,q}. 
\end{equation}
Consider the map
\begin{displaymath}
\begin{array}{rccc}
\psi: & \Phi & \rightarrow & \B_{n-2}\\
 & \beta&\mapsto  & (\alpha_{n-1,n}\beta)|_{\{1,\dots,n-2\}}.
\end{array}
\end{displaymath}
As in Proposition~\ref{prnewrec} one obtains that $\psi$ is a homomorphism
from $\Phi$ to $\B_{n-2}$, which maps $\alpha_{x,y}$ to the identity element
if $y\in\{n-1,n\}$ and to an element of corank $2$ otherwise. Furthermore, it is
easy to see that 
$\{\psi(\alpha_{x,y})\,:\,y<n-1\}$ satisfy all the assumptions of our statement
and hence by induction we obtain
\begin{displaymath}
\psi(\alpha_{i,j}\alpha_{s,t})=\psi(\alpha_{i,j})\psi(\alpha_{s,t})=
\psi(\alpha_{u,v})\psi(\alpha_{p,q})=\psi(\alpha_{u,v}\alpha_{p,q})
\end{displaymath}
if $j,v<n-1$. If $j\in\{n-1,n\}$ or $v\in\{n-1,n\}$ then \eqref{lbr1.1} is
straightforward since our ``parameters'' $\alpha_{x,y}^{(s)}$ do not affect
any part of \eqref{lbr1.1} at all. This implies \eqref{lbr1.1}.

Thus we are left to prove \eqref{lbr1.2} and \eqref{lbr1.3}.
If there exist $f,g\in M$ such that 
$(n-1)\equiv_{\alpha_{i,j}\alpha_{s,t}} f$ and
$n\equiv_{\alpha_{i,j}\alpha_{s,t}} g$, then we have
either $(n-1)\equiv_{\alpha_{i,j}\alpha_{s,t}}i$ and 
$n\equiv_{\alpha_{i,j}\alpha_{s,t}}u$ or $(n-1)\equiv_{\alpha_{i,j}\alpha_{s,t}}u$ and 
$n\equiv_{\alpha_{i,j}\alpha_{s,t}}i$. Combining this with \eqref{brbr1} we obtain either 
$(n-1)\equiv_{\alpha_{u,v}\alpha_{p,q}}i$ and $n\equiv_{\alpha_{u,v}\alpha_{p,q}}u$ or
$(n-1)\equiv_{\alpha_{u,v}\alpha_{p,q}}u$ and $n\equiv_{\alpha_{u,v}\alpha_{p,q}}i$, 
respectively.

Assume that there exist $f,g\in M$ such that 
$(n-1)\equiv_{\alpha_{i,j}\alpha_{s,t}} f'$ and
$n\equiv_{\alpha_{i,j}\alpha_{s,t}} g'$. Then we in fact have to prove that
\begin{equation}\label{brbr3}
\alpha_{i,j}^{(1)}\alpha_{s,t}^{(2)}=\alpha_{u,v}^{(1)}\alpha_{p,q}^{(2)}.
\end{equation}
If the pairs $(i,j)$, $(s,t)$, $(u,v)$, and $(p,q)$ either all belong or 
all do not belong to $\mathfrak{A}_n$, then our equality reduces to the 
obvious identities $x_1x_2=x_1x_2$ and $y_1y_2=y_1y_2$, respectively. 

Assume now that $(i,j)\in \mathfrak{A}_n$ and $(s,t)\not\in \mathfrak{A}_n$. If
$t<i$, then one obtains $(u,v)=(s,t)\not\in \mathfrak{A}_n$ and
$(p,q)=(i-2,j-2)\in \mathfrak{A}_n$.  If
$s\geq i$, then $(s,t)\not\in \mathfrak{A}_n$ means that $s$ and $t$ belong
to different sets from the $\{n,n-1\}$, $\{n-2,n-3\}$,\dots. From 
Corollary~\ref{c3}\eqref{c3.3} it follows that in this case $u$ and $v$ belong
to different sets from the $\{n,n-1\}$, $\{n-2,n-3\}$,\dots as well. This implies
that $(u,v)\not\in \mathfrak{A}_n$. It also follows that 
$(p,q)=(i,j)\in \mathfrak{A}_n$.  Finally, assume that $s<i\leq t$, then
$v-u>2$ and hence $(u,v)\not\in \mathfrak{A}_n$. Moreover, in the last case we also
have $(p,q)=(i,j)\in \mathfrak{A}_n$. The case $(i,j)\not\in \mathfrak{A}_n$ and 
$(s,t)\in \mathfrak{A}_n$ is analogous and we obtain that exactly one pair on the 
right hand side of \eqref{brbr2}  belongs to 
$\mathfrak{A}_n$. Then \eqref{brbr3} reduces to $x_1y_2=y_1x_2$, which is again 
the case. 

Finally, let us assume that there exist $f,g\in M$ such that 
$(n-1)\equiv_{\alpha_{i,j}\alpha_{s,t}} f$ and
$n\equiv_{\alpha_{i,j}\alpha_{s,t}} g'$. In this case we have either
$i=f$, $j=n-1$ and $\alpha_{i,j}$ is an idempotent or $u=f$, $v=n-1$, and 
$\alpha_{u,v}$ is an idempotent. Without loss of generality we assume that 
$\alpha_{i,j}$ is an idempotent. Now, using \eqref{brbr1},
we obtain $(n-1)\equiv_{\alpha_{u,v}\alpha_{p,q}} f$. Moreover, we obviously
have $\alpha_{u,v}=\alpha_{s,t}$ and $f\equiv_{\alpha_{s,t}}g'$. Observe that
$p=g$ and $q\in\{n-2,n-3\}$. If $q=n-3$, we have  $(n-2)\equiv_{\alpha_{p,q}}g'$.
If $q=n-2$, we have  $(n-3)\equiv_{\alpha_{p,q}}g'$. Hence in both cases we obtain
$n\equiv_{\alpha_{u,v}\alpha_{p,q}}g'$. For the case 
$(n-1)\equiv_{\alpha_{i,j}\alpha_{s,t}} g'$ and $n\equiv_{\alpha_{i,j}\alpha_{s,t}} f$
the arguments are similar. This completes the proof of Step~1.

{\em Step~2.} Now we go to elements of arbitrary coranks. Let
$\eta_1\cdot\dots\cdot\eta_k\in \Phi_k$. If $\zeta\in \B_n$
has corank $2$ then for every $\a\in \B_n$ we have
that $\mathrm{corank}(\a\zeta)$ is either $\mathrm{corank}(\a)$ or 
$\mathrm{corank}(\a)+2$. The same holds for $\mathrm{corank}(\zeta\a)$.
This implies that for every $i=1,\dots,k-1$ we have
$\mathrm{corank}(\eta_i\eta_{i+1})=4$. In particular, by Step~1 there exists
a unique pair, $(\eta'_i,\eta'_{i+1})\in\Gamma_1\times\Gamma_2$ such that
$(\eta'_i,\eta'_{i+1})\neq (\eta_i,\eta_{i+1})$ and
$\eta_i\eta_{i+1}=\eta'_i\eta'_{i+1}$, which, in particular, implies
that $\eta'_i\in\Gamma_i$ and $\eta'_{i+1}\in\Gamma_{i+1}$. This allows us to define the
involution $\mathfrak{i}^{(k)}_i$ on the set 
$\Gamma(k)=\Gamma_1\times\Gamma_2\times\dots\times\Gamma_k$ via 
\begin{displaymath}
\mathfrak{i}^{(k)}_i(\eta_1,\dots,\eta_k)=
(\eta_1,\dots,\eta_{i-1},\eta_i',\eta_{i+1}',\eta_{i+2},\dots,\eta_k).
\end{displaymath}
Note that for $(\eta_1,\dots,\eta_k),(\eta'_1,\dots,\eta'_k)\in \Gamma(k)$ such that
\begin{displaymath}
(\eta_1,\dots,\eta_k)=\mathfrak{i}^{(k)}_i(\eta'_1,\dots,\eta'_k)
\end{displaymath}
from Step~1 we have
\begin{equation}\label{eqinvol}
\eta_1\cdot\dots\cdot\eta_k=\eta'_1\cdot\dots\cdot\eta'_k.
\end{equation}

Fix $(\eta_1,\dots,\eta_k)\in\Gamma(k)$ and let $\a=\eta_1\cdot\dots\cdot\eta_k$. 
Let $X_1,\dots,X_k$ be the ordered collection of disjoint two-element subsets of 
$M$ such that for every $i=1,\dots,k$ the element $\eta_1\cdot\dots\cdot\eta_i$ 
contains $X_1,\dots,X_i$. Set
\begin{displaymath}
\Gamma(k,\alpha)=\{(\eta'_1,\dots,\eta'_k)\in\Gamma(k)\,:\,
\eta'_1\cdot\dots\cdot\eta'_k\mathcal{R}\alpha\}.
\end{displaymath}
We already know that $|\Gamma(k,\alpha)|=k!$, moreover, the arguments used to prove 
this allow us to define a bijective map, $\varphi:\Gamma(k,\alpha)\to\mathcal{S}_k$, 
in the following way: to an element, 
$(\eta'_1,\dots,\eta'_k)\in\Gamma(k,\alpha)$,
we associate $\sigma\in \S_k$ such that for every $i=1,\dots,k$ 
the element $\eta'_1\cdot\dots\cdot\eta'_i$ contains $X_{\sigma(1)},\dots,X_{\sigma(i)}$.
Let $\mathfrak{s}_i$ denote the simple transposition $(i,i+1)\in\S_k$. From the
definitions of $\varphi$ and $\mathfrak{i}^{(k)}_i$ for arbitrary $\beta\in 
\Gamma(k,\alpha)$ we have
\begin{equation}\label{eqbraid}
\varphi\circ \mathfrak{i}^{(k)}_i(\beta)=\varphi(\beta)\mathfrak{s}_i.
\end{equation}
For every $\sigma\in \S_k$ fix a reduced decomposition, 
$\mathfrak{s}_{a^{\sigma}_{l_{\sigma}}}\dots \mathfrak{s}_{a^{\sigma}_2}
\mathfrak{s}_{a^{\sigma}_{1}}$, 
of $\sigma$ and define
\begin{displaymath}
\a(\sigma)=\mathfrak{i}^{(k)}_{a^{\sigma}_1}\circ
\mathfrak{i}^{(k)}_{a^{\sigma}_2}\circ\dots\circ
\mathfrak{i}^{(k)}_{a^{\sigma}_{l_{\sigma}}}(\eta_1,\dots,\eta_k).
\end{displaymath}
From \eqref{eqbraid} we obtain $\varphi(\a(\sigma))=\sigma$, implying
$|\{\a(\sigma)\,:\,\sigma\in  \S_k\}|=k!$. This yields
$\{\a(\sigma)\,:\,\sigma\in  \S_k\}=\Gamma(k,\alpha)$ and hence
\eqref{eqinvol} implies that all elements of $\Gamma(k,\alpha)$ produce
the same element, namely $\a$, via multiplication. 
This means that $\overline{\Phi}$ contains precisely
one element in each $\mathcal{R}$--class and is exactly what 
we wanted to prove.

{\em Step~3.} Obviously, $\overline{\Phi}\subset\Phi$,
and we are left to show that $\overline{\Phi}=\Phi$. 
From the definition of $\overline{\Phi}$ it follows that
it is enough to show that $\a\eta\in \Phi_k$ 
for every $\a\in\Phi_k$ and $\eta\in\Gamma\setminus\Gamma_{k}$.
However, $\a\in\Phi_k$ implies that 
\begin{displaymath}
n'\equiv_{\tau}(n-1)',(n-2)'\equiv_{\tau}(n-3)',\dots,
(n-2(k-1))'\equiv_{\tau}(n-2(k-1)-1)',
\end{displaymath}
by the definition of $\Phi_k$. At the same time 
$\eta\in\Gamma\setminus\Gamma_{k}$ implies that $\eta=\alpha_{i,j}$ and $j>n-2k$.
Now Corollary~\ref{c3} and a direct calculation imply $\a\eta=\a$, which completes
the proof.
\end{proof}

\begin{theorem}\label{tmain2}
Let $\Lambda$ be a canonical $\mathcal{R}$--cross-section of $\B_n$. Then the 
generators $\alpha_{i,j}$, $1\leq i<j\leq n$, and the following relations:
\begin{enumerate}[(a)]
\item\label{tmain2.1} $\alpha_{s,t}\alpha_{i,j}=\alpha_{s,t}$ for all possible
$s,t,i$ and for all $j=n-1,n$,
\item\label{tmain2.2} $\alpha_{s,t}\alpha_{i,j}=\alpha_{u,v}\alpha_{x,y}(=\beta)$
for all $\beta\in\Lambda$, $\mathrm{corank}(\beta)=4$, with
appropriate $s,t,i,j,u,v,x,y$ given by Step~1 of the proof of  Proposition~\ref{p12},
\end{enumerate}
form a copresentation of the monoid $\Lambda$.
\end{theorem}

\begin{proof}
Let $\Psi$ denote the monoid, generated by 
$\gamma_{i,j}$, $1\leq i<j\leq n$, satisfying 
\begin{enumerate}[(a)]
\setcounter{enumi}{2}
\item\label{tmain2.1.1} 
$\gamma_{s,t}\gamma_{i,j}=\gamma_{s,t}$ for all possible
$s,t,i$ and for all $j=n-1,n$,
\item\label{tmain2.2.2} 
$\gamma_{s,t}\gamma_{i,j}=\gamma_{u,v}\gamma_{x,y}$
for all $s,t,i,j,u,v,x,y$ such that
$\alpha_{s,t}\alpha_{i,j}=\alpha_{u,v}\alpha_{x,y}$ is
an element of corank $4$.
\end{enumerate}
Let $\mathfrak{f}:\{\gamma_{i,j}:1\leq i<j\leq n\}\to 
\{\alpha_{i,j}:1\leq i<j\leq n\}$ be the bijection
given by $\mathfrak{f}(\alpha_{i,j})=\gamma_{i,j}$.
Then $\mathfrak{f}$ obviously extends to an epimorphism, 
$\overline{\mathfrak{f}}:\Psi\to \Lambda$.
To complete the proof it is enough to show that
$|\Psi|=|\Lambda|$. For $k\in\{1,\dots,\lfloor\frac{n}{2}\rfloor\}$ let
$\Psi^{(k)}$ denote the set $\{\gamma_{i,j}\,:\,j\leq n-2(k-1)\}$. Define
$\Psi_0=\{\mathrm{id}\}$ and for $k\in\{1,\dots,\lfloor\frac{n}{2}\rfloor\}$
let $\Psi_k$ denote the image of the multiplication map
\begin{displaymath}
\mathrm{mult}:\Psi^{(1)}\times\dots\times\Psi^{(k)}\to \Psi.
\end{displaymath}
Define $\overline{\Psi}=\cup_{k=0}^{\lfloor\frac{n}{2}\rfloor}\Psi_k$. Using
the relations from \eqref{tmain2.2.2} and the arguments analogous to those of
Proposition~\ref{decomposition} and Step~2 of Proposition~\ref{p12} one shows that
$|\Psi_k|=|\Lambda_k|$ for every $k\in\{1,\dots,\lfloor\frac{n}{2}\rfloor\}$.
Hence to complete the proof we have just to show that 
$\Psi=\overline{\Psi}$. To prove this it is enough to show that
for any $\gamma\in\Psi_k$ and $\eta\not\in \Psi^{(k+1)}$ we have
$\gamma\eta\in\Psi_k$. Assume that $\eta\in \Psi^{(l)}\setminus \Psi^{(l+1)}$ 
for some $l\leq k$, that is $\eta=\gamma_{u,v}$ and $v\in\{n-2(l-1),n-2(l-1)-1\}$.  Let $\gamma=\eta_1\cdot\dots\cdot\eta_k$ 
and $\eta_i\in\Psi^{(i)}$ for all $i$. 
Consider $\eta_1\cdot\dots\cdot\eta_k\eta$.
Using \eqref{tmain2.2.2} we have
$\eta_k\eta=\eta^{(1)}\eta'_k$, moreover,
$\mathfrak{f}(\eta_k)\mathfrak{f}(\eta)=
\mathfrak{f}(\eta^{(1)})\mathfrak{f}(\eta'_k)$.
Hence, using Corollary~\ref{c3}\eqref{c3.3}, we obtain
$\mathfrak{f}(\eta^{(1)})\in \Lambda_1^{(l-1)}\setminus\Lambda_1^{(l)}$ 
and $\mathfrak{f}(\eta'_k)\in \Lambda_1^{(k)}$,
which implies $\eta^{(1)}\in \Psi^{(l-1)}\setminus\Psi^{(l)}$ and 
$\eta'_k\in \Psi^{(k)}$. We proceed inductively and 
after $l-1$ steps we obtain
\begin{displaymath}
\eta_1\cdot\dots\cdot\eta_k\eta=
\eta_1\cdot\dots\cdot\eta_{k-l+1}\eta^{(l-1)}
\eta'_{k-l+2}\cdot\dots\cdot\eta'_k,
\end{displaymath}
where $\eta^{(l-1)}\in \Psi^{(1)}\setminus\Psi^{(2)}$ and
$\eta'_i\in\Psi^{(i)}$ for all $i>k-l+1$. Now 
\eqref{tmain2.1.1} implies $\eta_{k-l+1}\eta^{(l-1)}=\eta_{k-l+1}$
and hence $\eta_1\cdot\dots\cdot\eta_k\eta\in \Psi_k$. This completes
the proof.
\end{proof}

\section{The main results}\label{s4}

\begin{lemma}\label{l15}
Let $\Theta$ be an arbitrary $\mathcal{R}$--cross-section in $\B_n$.
Then $\Theta$ is generated (as a monoid) by elements of corank $2$.
\end{lemma}

\begin{proof}
Let $\Theta_1$ be the set of elements of corank $2$ in $\Theta$.
We have to show that $\alpha\in\langle\Theta_1\rangle$ for every
$\a\in\Theta$. We use the induction on $k$ such that 
$2k=\mathrm{corank}(\a)$, with the case $k=1$ being trivial. 
Let $\a\in\Theta$ be such that  $\mathrm{corank}(\a)=2k$, 
and let $i_1,\dots,i_k,j_1,\dots,j_k\in M$
be pairwise different and such that $i_s\equiv_{\a} j_s$ for all
$s=1,\dots,k$. Let $\beta$ be the element of $\Theta$ of
corank $2k-2$ such that $i_s\equiv_{\beta} j_s$ for all
$s=1,\dots,k-1$. By inductive assumption, $\b$ can be decomposed 
into a product of elements of corank $2$ from $\Theta$. Further,
since $\mathrm{corank}(\b)=2k-2$, we have that there exist
$u,v\in M$ such that $i_k\equiv_{\beta} u'$ and
$j_k\equiv_{\beta} v'$. Let further $\gamma$ be the element of
corank $2$ in $\Theta$ such that $u\equiv_{\gamma} v$. Then
$\beta\gamma\in\Theta$,  $\beta\gamma\mathcal{R}\alpha$ and thus 
$\beta\gamma=\a$ as $\Theta$ is an $\mathcal{R}$--cross-section.
This proves the induction step and completes the proof.
\end{proof}

The group $\S_n$ acts by isomorphisms on $\B_n$ via conjugation. 
This action induces an action of $\S_n$ on the
set of all $\mathcal{R}$--cross-sections of $\B_n$. Our next step towards
the classification of all $\mathcal{R}$--cross-sections of $\B_n$ is
the following statement:

\begin{proposition}\label{p16}
Every $\mathcal{R}$--cross-section of $\B_n$ is $\S_n$-conjugated to a
canonical $\mathcal{R}$--cross-section.
\end{proposition}

\begin{proof}
Let $\Theta$ be an $\mathcal{R}$--cross-section of $\B_n$.
To prove our statement it is certainly enough to show that there exists a
sequence, $X_1,\dots,X_{\lfloor\frac{n}{2}\rfloor}$, of disjoint
two-element subsets of $M$ such that for every 
$k\in\{1,\dots,\lfloor\frac{n}{2}\rfloor\}$ each element $\a\in\Theta$ of 
corank $2k$ contains $X'_1,X'_2,\dots,X'_k$. 

To prove our statement we will need two lemmas.

\begin{lemma}\label{lh1}
All elements of $\Theta$ of corank $2$ belong to the same
$\mathcal{L}$-class.
\end{lemma}

\begin{proof}
Let $\a\in\Theta$ be an element of corank $2$
and $i'\equiv_{\a}j'$ for some $i,j\in M$.  Consider the element $\beta\in\Theta$ 
of corank $2$ such that $i\equiv_{\b}j$. Then $\a\beta$ has corank $2$, moreover,
$\a\beta\mathcal{R}\a$. Thus $\a\beta=\a$, which implies that $\b$ must coincide
with the (unique) idempotent $\pi_{\{i,j\}}$ of corank $2$, which satisfies 
$s\equiv_{\pi_{\{i,j\}}}s'$ for all $s\not\in\{i,j\}$.

Now assume that there are $\a,\b\in\Theta$ such that
$\mathrm{corank}(\a)=\mathrm{corank}(\b)=2$ and the elements  $\a$ and $\b$ 
belong to different $\mathcal{L}$--classes. We have to consider two cases. 
First we assume that there exist different $i,j,t\in M$ such that $i'\equiv_{\a}j'$ 
and  $j'\equiv_{\b}t'$. In this case the above observation implies 
$\pi_{\{i,j\}},\pi_{\{j,t\}}\in\Theta$. But then 
$\pi_{\{i,j\}}\pi_{\{j,t\}}\mathcal{R}\pi_{\{i,j\}}$ and
$\pi_{\{i,j\}}\pi_{\{j,t\}}\neq\pi_{\{i,j\}}$
(see Figure~\ref{fig:f5}), which contradicts the fact that 
\begin{figure}
\special{em:linewidth 0.4pt}
\unitlength 0.80mm
\linethickness{0.4pt}
\begin{picture}(155.00,35.00)
\drawline(30.00,10.00)(30.00,30.00)
\drawline(30.00,10.00)(60.00,10.00)
\drawline(30.00,30.00)(60.00,30.00)
\drawline(60.00,10.00)(60.00,30.00)
\drawline(25.00,15.00)(30.00,15.00)
\drawline(25.00,20.00)(30.00,20.00)
\drawline(25.00,25.00)(30.00,25.00)
\drawline(65.00,15.00)(60.00,15.00)
\drawline(65.00,20.00)(60.00,20.00)
\drawline(65.00,25.00)(60.00,25.00)
\drawline(30.00,25.00)(60.00,25.00)
\drawline(30.00,15.00)(32.00,15.00)
\drawline(30.00,20.00)(32.00,20.00)
\drawline(32.00,20.00)(32.00,15.00)
\drawline(60.00,15.00)(58.00,15.00)
\drawline(60.00,20.00)(58.00,20.00)
\drawline(58.00,15.00)(58.00,20.00)
\drawline(70.00,10.00)(70.00,30.00)
\drawline(70.00,10.00)(100.00,10.00)
\drawline(70.00,30.00)(100.00,30.00)
\drawline(100.00,10.00)(100.00,30.00)
\drawline(65.00,15.00)(70.00,15.00)
\drawline(65.00,20.00)(70.00,20.00)
\drawline(65.00,25.00)(70.00,25.00)
\drawline(105.00,15.00)(100.00,15.00)
\drawline(105.00,20.00)(100.00,20.00)
\drawline(105.00,25.00)(100.00,25.00)
\drawline(70.00,15.00)(100.00,15.00)
\drawline(70.00,20.00)(72.00,20.00)
\drawline(70.00,25.00)(72.00,25.00)
\drawline(72.00,25.00)(72.00,20.00)
\drawline(100.00,20.00)(98.00,20.00)
\drawline(100.00,25.00)(98.00,25.00)
\drawline(98.00,20.00)(98.00,25.00)
\drawline(115.00,10.00)(115.00,30.00)
\drawline(115.00,10.00)(145.00,10.00)
\drawline(115.00,30.00)(145.00,30.00)
\drawline(145.00,10.00)(145.00,30.00)
\drawline(110.00,15.00)(115.00,15.00)
\drawline(110.00,20.00)(115.00,20.00)
\drawline(110.00,25.00)(115.00,25.00)
\drawline(150.00,15.00)(145.00,15.00)
\drawline(150.00,20.00)(145.00,20.00)
\drawline(150.00,25.00)(145.00,25.00)
\drawline(115.00,25.00)(145.00,15.00)
\drawline(115.00,15.00)(117.00,15.00)
\drawline(115.00,20.00)(117.00,20.00)
\drawline(117.00,15.00)(117.00,20.00)
\drawline(145.00,20.00)(143.00,20.00)
\drawline(145.00,25.00)(143.00,25.00)
\drawline(143.00,20.00)(143.00,25.00)
\put(107.50,20.00){\makebox(0,0)[cc]{$=$}}
\put(45.00,5.00){\makebox(0,0)[cc]{$\pi_{i,j}$}}
\put(85.00,5.00){\makebox(0,0)[cc]{$\pi_{j,t}$}}
\put(130.00,5.00){\makebox(0,0)[cc]{$\pi_{i,j}\pi_{j,t}$}}
\end{picture}
\caption{$\pi_{\{i,j\}}\pi_{\{j,t\}}\mathcal{R}\pi_{\{i,j\}}$ and
$\pi_{\{i,j\}}\pi_{\{j,t\}}\neq\pi_{\{i,j\}}$.}\label{fig:f5}
\end{figure}
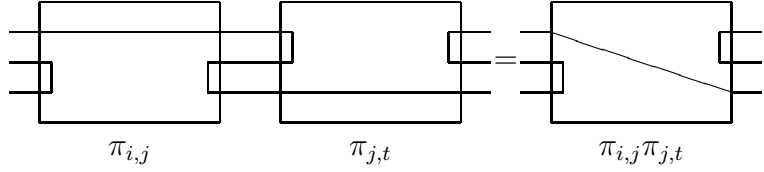
$\Theta$ is an $\mathcal{R}$--cross-section. In the second case we assume that 
there exist pairwise different $i,j,s,t\in M$ such that $i'\equiv_{\a}j'$ and 
$s'\equiv_{\b}t'$. In this case the above observation implies 
$\pi_{\{i,j\}},\pi_{\{s,t\}}\in\Theta$. Let $\gamma\in\Theta$ be such that
$j\equiv_{\gamma}s$. Then it is easy to see that 
$\pi_{\{i,j\}}\gamma\mathcal{R}\pi_{\{i,j\}}$ and thus $i'\equiv_{\gamma}j'$
since $\Theta$ is an $\mathcal{R}$--cross-section. Analogously we obtain
$s'\equiv_{\gamma}t'$, which contradicts the fact that $\mathrm{corank}(\gamma)=2$.
This completes the proof. 
\end{proof}

Lemma~\ref{lh1} and Theorem~\ref{tgreen}\eqref{tgreen.1} say 
that there exists a two-element subset, $X_1$, of $M$
such that $X'_1$ belongs to every element of $\Theta$ of corank $2$.
Since the elements of corank $2$ generate $\Theta$ by Lemma~\ref{l15} we even
have that $X'_1$ belongs to every element of $\Theta$ of corank at least $2$.
Let $\pi$ denote the unique element of $\Theta$ of corank $2$, containing $X_1$. It
is obviously an idempotent and $i\equiv_{\pi}i'$ for all $i\not\in X_1$. 
Define $\overline{\varphi}:\Theta\to\Theta$ via $\overline{\varphi}(\a)=\pi\a\pi$,
$\a\in\Theta$, and set $N=\overline{\varphi}(\Theta)$.

\begin{lemma}\label{lh2}
\begin{enumerate}[(i)]
\item\label{lh2.1} $\overline{\varphi}(\a)=\pi\a$ for all $\a\in\Theta$.
\item\label{lh2.2} The set $M\setminus X_1$ is $\overline{\varphi}(\a)$-invariant
for all $\a\in\Theta$, in particular, the element
$\varphi(\a)=(\overline{\varphi}(\a))|_{M\setminus X_1}$ is a well-defined
element of $\B(M\setminus X_1)$.
\item\label{lh2.3} The map $\varphi:\Theta\to \B(M\setminus X_1)$ is
a homomorphism and its image is an $\mathcal{R}$--cross-section of 
$\B(M\setminus X_1)$.
\end{enumerate}
\end{lemma}

\begin{proof}
Analogous to that of Proposition~\ref{prnewrec} and hence is left to the reader.
\end{proof}

After Lemma~\ref{lh2} we can use induction on $n$. First we observe that for
$n=1,2$ our statement is trivial. Let $n>2$. Since 
$\varphi(\Theta)$ is an $\mathcal{R}$--cross-section of 
$\B(M\setminus X_1)$ by Lemma~\ref{lh2}\eqref{lh2.3}, by induction there
exists a sequence, $X_2,\dots,X_{\lfloor\frac{n}{2}\rfloor}$, or disjoint
two-element subsets of $M\setminus X_1$ such that for every 
$k\in\{1,\dots,\lfloor\frac{n}{2}\rfloor-1\}$ each element $\beta\in\varphi(\Theta)$ 
of corank $2k$ contains $X'_2,\dots,X'_{k+1}$. 

Now let $k\in\{1,\dots,\lfloor\frac{n}{2}\rfloor\}$
and $\a\in \Theta$ be an element of corank $2k$. Let
$X_1=\{x,y\}$. Assume first that either
$x$ or $y$ (or both) does not belong to a line in $\a$. 
Then it is easy to see that
the element $\overline{\varphi}(\a)$ also has corank $2k$ and thus
$\varphi(\a)$ has corank $2(k-1)$. Thus $\varphi(\a)$ contains
$X'_2,\dots,X'_{k}$, which implies that $\a$ contains
$X'_1,\dots,X'_{k}$.

If both, $x$ and $y$, belong to some lines in $\a$, then 
the element $\overline{\varphi}(\a)$ has corank 
$2k+2$ and thus $\varphi(\a)$ has corank $2k$. Therefore,
by the induction hypothesis,
$\varphi(\a)$ contains $X'_2,\dots,X'_{k+1}$ and thus 
$\overline{\varphi}(\a)$ contains $X'_1,\dots,X'_{k+1}$.
Using Lemma~\ref{l15} we can write $\a=\eta\beta$, where
$\eta$ has corank $2$ and $\beta$ has corank $2k-2$. Using induction on 
$k$ (the case $k=1$ is given by Lemma~\ref{lh1}), we can assume that
$\beta$ contains $X'_1,\dots,X'_{k-1}$ and hence so does $\a$. Now 
we have to find out if $\a$ contains $X'_{k}$ or $X'_{k+1}$ (it must contain
exactly one of them since $\mathrm{corank}(\a)=2k$). Assume that
$\a$ contains $X'_{k+1}$ and take any $\beta$ of corank $2k$, containing
$\{x,y\}$. Then $\beta$ must contain $X'_1,\dots,X'_{k}$ by the previous 
paragraph. Let $X_k=\{a,b\}$ and $X_{k+1}=\{u,v\}$ and note that
$\{a,b\}\cap \{u,v\}=\varnothing$ by our assumptions. 
Consider the unique element $\gamma$ of $\Theta$ of corank $2k$, which 
contains $X_1,\dots,X_{k-1}$ and $\{a,u\}$. Then $\gamma$ has 
to contain $X'_1, \dots, X'_k$. It is easily checked that
$\alpha\gamma\mathcal{R}\alpha$ and hence 
$\alpha\gamma=\alpha$, a contradiction.
Hence $\alpha$ must contain $X'_k$ and the proof is complete.
\end{proof}

Recall that for  $i\neq j\in M$ we denote by $\mathfrak{s}_{i}$ the transposition
$(i,i+1)\in\S_n\subset\B_n$. A canonical $\mathcal{R}$--cross-section,
$\Lambda$ will be called {\em regular} provided that
$\a_{i,j}^{(l)}=\mathrm{id}$ for all possible $i,j,l$.
A canonical $\mathcal{R}$--cross-section,
$\Lambda$ will be called {\em alternating} provided that 
$\a_{i,j}^{(l)}=(1,2)$ for all possible $(i,j)\not\in\mathfrak{A}_n$ 
and all possible $l$, and $\a_{i,i+1}^{(l)}=\mathrm{id}$ for all 
$(i,i+1)\in\mathfrak{A}_n$ and for all possible $l$. Now we are ready to
formulate our main result.

\begin{theorem}\label{tmain1}
The stabilizer in $\S_n$ of any canonical $\mathcal{R}$--cross-section consists of
$\mathrm{id}\in\S_n$ and
$\mathfrak{s}_{n-1,n}\prod_{(i,i+1)\in\mathfrak{A}_n}\mathfrak{s}_{i,i+1}\in\S_n$.
\begin{enumerate}[(i)]
\item\label{tmain1.1} For $n=1$ we have one trivial $\mathcal{R}$--cross-section.
\item\label{tmain1.2} For $n=2$ we have only one $\mathcal{R}$--cross-section, moreover,
it is canonical and consists of $\alpha_{1,2}$ and $\mathrm{id}$.
\item\label{tmain1.3} For $n=3$ we have $3$ different $\mathcal{R}$--cross-sections, namely
$\Lambda$, $\mathfrak{s}_{1,2}\Lambda\mathfrak{s}_{1,2}$, and 
$\mathfrak{s}_{1,3}\Lambda\mathfrak{s}_{1,3}$, where $\Lambda$ is the canonical 
cross-section described in Proposition~\ref{p10}\eqref{p10.3}. 
\item\label{tmain1.4} For $n=4$ we have $12$ different $\mathcal{R}$--cross-sections,
each of which is $\S_4$-conjugated to the regular canonical
$\mathcal{R}$--cross-section. 
\item\label{tmain1.5} For $n=5$ we have $2\cdot 5!$ different $\mathcal{R}$--cross-sections,
each of which is $\S_5$-conjugated to a canonical $\mathcal{R}$--cross-section
from Proposition~\ref{p10}\eqref{p10.5} satisfying $\alpha_{1,2}^{(1)}=\mathrm{id}$.
\item\label{tmain1.6} For $n=6$ we have $2\cdot 6!$ different $\mathcal{R}$--cross-section, 
each of which is $\S_6$-conjugated to a canonical $\mathcal{R}$--cross-section
from Proposition~\ref{p10}\eqref{p10.6} satisfying $\alpha_{1,2}^{(1)}=\mathrm{id}$
and $\alpha_{1,2}^{(2)}=\mathrm{id}$. 
\item\label{tmain1.7} For $n\geq 7$ we have $n!$ different $\mathcal{R}$--cross-section.
Half of them are $\S_n$-conjugated to the regular canonical
$\mathcal{R}$--cross-section, and half of them are  $\S_n$-con\-ju\-gat\-ed
to the alternating canonical $\mathcal{R}$--cross-section. 
\end{enumerate}
\end{theorem}

\begin{proof}
Since every $\mathcal{R}$--cross-section is $\S_n$-conjugated to a canonical
$\mathcal{R}$--cross-section by Proposition~\ref{p16}, to complete the classification
we have to determine which canonical $\mathcal{R}$--cross-sections of $\B_n$
are $\S_n$-conjugated to each other. Let $\Lambda$ be a canonical $\mathcal{R}$--cross-section 
and $\sigma\in\S_n\subset \B_n$. Assume that $\sigma$ stabilizes $\Lambda$. In particular,
$\sigma^{-1}\Lambda\sigma$ is a canonical $\mathcal{R}$--cross-section, which
forces $\sigma(n)\in\{n-1,n\}$, $\sigma(n-1)\in\{n-1,n\}$, $\sigma(n-2)\in\{n-3,n-2\}$,
and so on. This implies that $\sigma$ must stabilize either $\alpha_{1,2}$
(if $n$ is even) or $\alpha_{2,3}$ (if $n$ is odd). A direct calculation shows that 
this is possible if and only if $\sigma=\mathrm{id}$ or
$\sigma=\mathfrak{s}_{n-1}\prod_{(i,i+1)\in\mathfrak{A}_n}\mathfrak{s}_{i}$.  

Consider the case $n\geq 7$. From Proposition~\ref{p11} and Proposition~\ref{p12} we 
have that there are exactly $2^{\lfloor\frac{n}{2}\rfloor}$ canonical 
$\mathcal{R}$--cross-sections. The arguments from the previous paragraph imply that 
the number of canonical conjugates of a given canonical 
$\mathcal{R}$--cross-section is exactly 
$2^{\lfloor\frac{n}{2}\rfloor-1}$, and the $\S_n$-orbit of every canonical 
$\mathcal{R}$--cross-section has size $\frac{n!}{2}$. Using the parity arguments 
it is also easy to show that the regular and the alternating canonical 
$\mathcal{R}$--cross-sections are not $\S_n$-conjugated (in Theorem~\ref{tmain3}
below we even prove that they are not isomorphic). This completes the proof
in the case $n\geq 7$. The case $n\leq 6$ can be treated using analogous 
arguments and Proposition~\ref{p10}. This is left to the reader.
\end{proof}

We remark once more that the classification of $\mathcal{L}$--cross-sections
in $\B_n$ is obtained by applying ${}^*$ to the statement of 
Theorem~\ref{tmain1}.

\begin{theorem}\label{tmain3}
For all $n\geq 6$ the regular $\mathcal{R}$--cross-section $\Lambda$ and the
alternating $\mathcal{R}$--cross-section $\Gamma$ in $\B_n$ are not isomorphic 
as monoids.
\end{theorem}

To prove Theorem~\ref{tmain3} we will need the following lemma

\begin{lemma}\label{l31}
Let $\varphi:\Lambda\to\Gamma$ be an isomorphism. Then for every
$\a\in\Lambda$ we have $\mathrm{rank}(\a)=\mathrm{rank}(\varphi(\a))$ and
$\mathrm{strank}(\a)=\mathrm{strank}(\varphi(\a))$.
\end{lemma}

\begin{proof}
For $k=0,1,\dots,\lfloor\frac{n}{2}\rfloor$ denote by
$\Lambda_k$ and $\Gamma_k$ the set of all element of corank
$2k$ in $\Lambda$ and $\Gamma$, respectively. Obviously,
$\phi:\Lambda_0\to\Gamma_0$. Now we claim
that $\Lambda_1$ is the unique irreducible system of generators 
for $\Lambda$ (as a monoid), and  $\Gamma_1$ is the unique irreducible 
system of generators for $\Gamma$. We prove the statement for $\Lambda$ and 
for $\Gamma$ the arguments are the same. That $\Lambda_1$ generates
$\Lambda$ we know from Proposition~\ref{prnew1}\eqref{prnew1.3}. If
$\eta_1,\eta_2\in \Lambda_1$ are such that $\eta_1\eta_2\in\Lambda_1$
then Proposition~\ref{prnew1}\eqref{prnew1.1} implies $\eta_1\eta_2=\eta_1$.
Hence $\Lambda_1$ is irreducible. On the other hand, if $N$ is an
irreducible system of generators for $\Lambda$, then $N$ generates,
in particular, $\Lambda_1$, which means that $N\cap \Lambda_1$ generates
$\Lambda_1$. By Proposition~\ref{prnew1}\eqref{prnew1.3} we have that
$N\cap \Lambda_1$ generates $\Lambda$ and hence coincides with $N$ by
the minimality of $N$. The minimality of $\Lambda_1$ now implies $N=\Lambda_1$.

From the above arguments we deduce that $\phi:\Lambda_1\to\Gamma_1$ is
a bijection (since both sets have the same cardinality). From
Proposition~\ref{prnew1}\eqref{prnew1.2} it follows that the elements of
$\Lambda_k$ are characterized as those elements in 
$\Lambda$, whose shortest possible decomposition into a product of
elements from $\Lambda_1$ has length $k$. Since $\phi:\Lambda_1\to\Gamma_1$,
it  follows that $\phi:\Lambda_k\to\Gamma_k$ for all $k$ and hence $\phi$ 
preserves the coranks and thus the ranks of the elements.

Let $\a^i=\pi_{\a}\in\Lambda$ be an idempotent. Then 
$\varphi(\pi_{\a})=\pi_{\varphi(\a)}$ and thus 
$\mathrm{strank}(\a)=\mathrm{strank}(\varphi(\a))$ follows from the equality
$\mathrm{rank}(\pi_{\a})=\mathrm{rank}(\varphi(\pi_{\a}))$ proved above.
\end{proof}

\begin{proof}[Proof of Theorem~\ref{tmain3}.]
Suppose that  $\varphi:\Lambda\to\Gamma$ is an isomorphism of the regular
cross-section $\Lambda$ onto the alternating cross-section $\Gamma$.
We will write $\alpha_{i,j}$ for the elements from $\Lambda$ and
$\tilde{\alpha}_{i,j}$ for the corresponding elements from $\Gamma$.
In particular, $\tilde{\alpha}_{i,j}=\alpha_{i,j}$ for all $(i,j)\in\mathfrak{A}_n$
by the definition of $\Gamma$.

Let $n\geq 6$ be even. Then $\alpha_{1,2}$ (resp. $\tilde{\alpha}_{1,2}$)
is the unique element  in $\Lambda$ (resp. $\Gamma$) satisfying 
$\mathrm{rank}(\alpha_{1,2})=n-2$ and $\mathrm{strank}(\alpha_{1,2})=0$. 
From Lemma~\ref{l31} it follows that $\varphi(\a_{1,2})=\tilde{\a}_{1,2}$. 
A direct calculation shows that in $\Lambda$ we have the relation
$\alpha_{1,2}\alpha_{2,3}=\alpha_{4,5}\alpha_{1,2}$. Applying $\varphi$ gives
\begin{equation}\label{t3e2}
\tilde{\alpha}_{1,2}\varphi(\alpha_{2,3})=\varphi(\alpha_{4,5})\tilde{\alpha}_{1,2}.
\end{equation}
There are exactly $5$ elements in $\Gamma$ which have rank
$n-2$ and stable rank $2$, namely $\tilde{\alpha}_{2,3}$, $\tilde{\alpha}_{1,4}$,
$\tilde{\alpha}_{1,3}$, $\tilde{\alpha}_{2,4}$, and 
$\tilde{\alpha}_{3,4}$. Hence, by Lemma~\ref{l31},
$\varphi(\alpha_{2,3})$ must coincide with one of these elements. We consider
all these cases separately.

{\em Case 1: $\varphi(\alpha_{2,3})=\tilde{\alpha}_{2,3}$.}
Then \eqref{t3e2} implies $\varphi(\alpha_{4,5})=\tilde{\alpha}_{4,5}$ by a direct calculation.
Further, we also have the following relation in $\Lambda$:
\begin{equation}\label{t3e3}
\alpha_{2,3}\alpha_{2,3}=\alpha_{4,5}\alpha_{2,3},
\end{equation}
which implies $\tilde{\alpha}_{2,3}\tilde{\alpha}_{2,3}=\varphi(\alpha_{4,5})\tilde{\alpha}_{2,3}$.
However, a direct calculation shows that in this case we have
$\varphi(\alpha_{4,5})=\tilde{\alpha}_{4,6}$, a contradiction.

{\em Case 2: $\varphi(\alpha_{2,3})=\tilde{\alpha}_{1,4}$.}
Then \eqref{t3e2} implies $\varphi(\alpha_{4,5})=\tilde{\alpha}_{3,6}$, and
\eqref{t3e3} implies $\varphi(\alpha_{4,5})=\tilde{\alpha}_{3,5}$ by a direct calculation.
A contradiction.

{\em Case 3: $\varphi(\alpha_{2,3})=\tilde{\alpha}_{1,3}$.}
Then \eqref{t3e2} implies $\varphi(\alpha_{4,5})=\tilde{\alpha}_{3,5}$, and
\eqref{t3e3} implies $\varphi(\alpha_{4,5})=\tilde{\alpha}_{4,6}$ by a direct calculation.
A contradiction.

{\em Case 4: $\varphi(\alpha_{2,3})=\tilde{\alpha}_{2,4}$.}
Then \eqref{t3e2} implies $\varphi(\alpha_{4,5})=\tilde{\alpha}_{4,6}$, and
\eqref{t3e3} implies $\varphi(\alpha_{4,5})=\tilde{\alpha}_{3,5}$ by a direct calculation.
A contradiction.

{\em Case 5: $\varphi(\alpha_{2,3})=\tilde{\alpha}_{3,4}$.}
A direct calculation gives us $\alpha_{2,3}\alpha_{1,2}=\alpha_{1,4}\alpha_{1,2}$.
Applying $\varphi$ we have $\tilde{\alpha}_{3,4}\tilde{\alpha}_{1,2}=
\varphi(\alpha_{1,4})\tilde{\alpha}_{1,2}$. This implies
$\varphi(\alpha_{1,4})=\tilde{\alpha}_{1,2}=\varphi(\alpha_{1,2})$, which is not possible since
$\varphi$ is bijective. 

Hence $\varphi$ can not exist for even $n\geq 6$.

Now let $n\geq 7$ be odd. In this case we have exactly $3$ elements in $\Lambda$
(resp. $\Gamma$) of rank $n-2$ and stable rank $1$, namely $\alpha_{1,2}$,
$\alpha_{2,3}$ and $\alpha_{1,3}$ (resp. $\tilde{\alpha}_{1,2}$,
$\tilde{\alpha}_{2,3}$ and $\tilde{\alpha}_{1,3}$). Using Lemma~\ref{l31}, we
obtain that one of the following $6$ cases must occur.

{\em Case 1: $\varphi(\alpha_{1,2})=\tilde{\alpha}_{1,2}$,
$\varphi(\alpha_{1,3})=\tilde{\alpha}_{1,3}$ and 
$\varphi(\alpha_{2,3})=\tilde{\alpha}_{2,3}$.} 
We have the following relations in $\Lambda$:
\begin{gather}
\alpha_{1,2}\alpha_{1,2}=\alpha_{3,4}\alpha_{1,2}\label{t3e6},\\
\alpha_{3,4}\alpha_{3,4}=\alpha_{5,6}\alpha_{3,4}\label{t3e7},\\
\alpha_{2,3}\alpha_{3,4}=\alpha_{5,6}\alpha_{2,3}\label{t3e8}.
\end{gather}
Applying $\varphi$ and using a direct calculation, we obtain
$\varphi(\a_{3,4})=\tilde{\alpha}_{3,5}$ from \eqref{t3e6},
$\varphi(\a_{5,6})=\tilde{\alpha}_{4,6}$ from \eqref{t3e7}, and
$\varphi(\a_{5,6})=\tilde{\alpha}_{5,7}$ from \eqref{t3e8}, a contradiction.

{\em Case 2: $\varphi(\alpha_{1,2})=\tilde{\alpha}_{1,3}$,
$\varphi(\alpha_{1,3})=\tilde{\alpha}_{1,2}$ and 
$\varphi(\alpha_{2,3})=\tilde{\alpha}_{2,3}$.} 
Applying $\varphi$ and using a direct calculation, we obtain
$\varphi(\a_{3,4})=\tilde{\alpha}_{2,4}$ from \eqref{t3e6},
$\varphi(\a_{5,6})=\tilde{\alpha}_{5,7}$ from \eqref{t3e7}, and
$\varphi(\a_{5,6})=\tilde{\alpha}_{4,6}$ from \eqref{t3e8}, a contradiction.

{\em Case 3: $\varphi(\alpha_{1,2})=\tilde{\alpha}_{2,3}$,
$\varphi(\alpha_{1,3})=\tilde{\alpha}_{1,3}$ and 
$\varphi(\alpha_{2,3})=\tilde{\alpha}_{1,2}$.}
We have the following relations in $\Lambda$:
\begin{gather}
\alpha_{2,3}\alpha_{2,3}=\alpha_{4,5}\alpha_{2,3}\label{t3e9},\\
\alpha_{1,2}\alpha_{2,3}=\alpha_{4,5}\alpha_{1,2}\label{t3e10}.
\end{gather}
Applying $\varphi$ and using a direct calculation, we obtain
$\varphi(\a_{4,5})=\tilde{\alpha}_{3,5}$ from \eqref{t3e9} and
$\varphi(\a_{4,5})=\tilde{\alpha}_{1,4}$ from \eqref{t3e10}, a contradiction.

{\em Case 4: $\varphi(\alpha_{1,2})=\tilde{\alpha}_{1,3}$,
$\varphi(\alpha_{1,3})=\tilde{\alpha}_{2,3}$ and 
$\varphi(\alpha_{2,3})=\tilde{\alpha}_{1,2}$.}
Applying $\varphi$ and using a direct calculation, we obtain
$\varphi(\a_{4,5})=\tilde{\alpha}_{3,5}$ from \eqref{t3e9} and
$\varphi(\a_{4,5})=\tilde{\alpha}_{2,5}$ from \eqref{t3e10}, a contradiction.

{\em Case 5: $\varphi(\alpha_{1,2})=\tilde{\alpha}_{2,3}$,
$\varphi(\alpha_{1,3})=\tilde{\alpha}_{1,2}$ and 
$\varphi(\alpha_{2,3})=\tilde{\alpha}_{1,3}$.}
Applying $\varphi$ and using a direct calculation, we obtain
$\varphi(\a_{4,5})=\tilde{\alpha}_{2,4}$ from \eqref{t3e9} and
$\varphi(\a_{4,5})=\tilde{\alpha}_{1,5}$ from \eqref{t3e10}, a contradiction.

{\em Case 6: $\varphi(\alpha_{1,2})=\tilde{\alpha}_{1,2}$,
$\varphi(\alpha_{1,3})=\tilde{\alpha}_{2,3}$ and 
$\varphi(\alpha_{2,3})=\tilde{\alpha}_{1,3}$.}
Applying $\varphi$ and using a direct calculation, we obtain
$\varphi(\a_{4,5})=\tilde{\alpha}_{2,4}$ from \eqref{t3e9} and
$\varphi(\a_{4,5})=\tilde{\alpha}_{3,4}$ from \eqref{t3e10}, a contradiction.

Hence $\varphi$ can not exist for odd $n\geq 7$ either. This completes the proof.
\end{proof}

We remark that for $n\leq 4$ all $\mathcal{R}$--cross-sections in $\B_n$ are conjugated and
hence isomorphic. For $n=5$ the isomorphism of all canonical (and hence of all)
$\mathcal{R}$--cross-sections in $\B_5$ can be shown by a direct calculation.

\section{On $\mathcal{D}$-- and $\mathcal{H}$--cross-sections in $\B_n$}\label{s5}

It is of course a natural question what one can say about the 
$\mathcal{D}$-- and $\mathcal{H}$--cross-sections in $\B_n$. For
$\mathcal{H}$--cross-sections the answer is very easy.

\begin{proposition}\label{ps5.1}
For $n=1,2,3$ the semigroup $\B_n$ contains a unique 
$\mathcal{H}$--cross-section. This cross-section consists of all idempotents
in $\B_n$. For $n\geq 4$ the semigroup $\B_n$ does not contain any 
$\mathcal{H}$--cross-section.
\end{proposition}

\begin{proof}
For $n=1,2,3$ the statement is easily checked. For $n\geq 4$ we first observe
that an $\mathcal{H}$--cross-section must contain all idempotents of
the semigroup. A direct calculation shows that for $n=4$ the semigroup,
generated by all idempotents, is not an $\mathcal{H}$--cross-section. Moreover,
it contains $\mathcal{H}$--classes with more than one element, which shows
that no $\mathcal{H}$--cross-sections exist.
Using the canonical embedding $\B_4\hookrightarrow \B_n$, this also 
implies the same statement for all $n>4$.
\end{proof}

On the other hand, for $\mathcal{D}$--cross-sections we will now show that the
problem of their classification contains, as a sub-problem, the problem
of classification of all $\mathcal{D}$--cross-sections in the symmetric 
inverse semigroup $\mathcal{IS}_m$, where $m=\lfloor\frac{n}{2}\rfloor$. The
latter problem is still open, see \cite{GM3}.

\begin{proposition}\label{ps5.2}
Let $m=\lfloor\frac{n}{2}\rfloor$. Let $\Gamma$ be
a $\mathcal{D}$--cross-section in $\mathcal{IS}_m$. For
$f\in \Gamma$ define $\alpha_f\in\B_n$ as follows:
\begin{enumerate}
\item $(n-2i+1)\equiv_{\alpha_f}(n-2f(i)+1)'$ if $i\in\mathrm{dom}(f)$;
\item $(n-2i+2)\equiv_{\alpha_f}(n-2f(i)+2)'$ if $i\in\mathrm{dom}(f)$;
\item $(n-2i+1)\equiv_{\alpha_f}(n-2i+2)$ if $i\not\in\mathrm{dom}(f)$;
\item $(n-2i+1)'\equiv_{\alpha_f}(n-2i+2)'$ if $i\not\in\mathrm{ran}(f)$;
\item $1\equiv_{\alpha_f}1'$ if $n$ is odd.
\end{enumerate}
Then $\{\alpha_f\,:\,f\in\Gamma\}$ is a $\mathcal{D}$--cross-section in $\B_n$.
\end{proposition}

\begin{proof}
From the definition of $\alpha_f$ by a direct calculation it follows that  
$f\mapsto \alpha_f$ is a homomorphism, which implies that 
$\{\alpha_f\,:\,f\in\Gamma\}$ is a semigroup. Using 
Theorem~\ref{tgreen}\eqref{tgreen.4} one easily checks that 
$\{\alpha_f\,:\,f\in\Gamma\}$ is a $\mathcal{D}$--cross-section. The
statement follows.
\end{proof}

\vspace{0.2cm}

\begin{center}
\bf Acknowledgments
\end{center}
\vspace{0.2cm}

The paper was written during the visit of the first author to Uppsala University, 
which was supported by the Swedish Institute. The financial support of the  Swedish 
Institute and the hospitality of Uppsala University are gratefully acknowledged. 
For the third author the research was partially supported by the Swedish Research 
Council. We thank the referee for pointing out several inaccuracies in the original
version of the paper and for many useful comments which led to the improvements in the 
paper.

\vspace{0.3cm}

\noindent
G.K.: Algebra, Department of Mathematics and Mechanics, Kyiv Taras
Shevchenko University, 64 Volodymyrska st., 01033 Kyiv, UKRAINE,
e-mail: {\tt akudr\symbol{64}univ.kiev.ua}
\vspace{0.3cm}

\noindent
V.Mal.: Algebra, Department of Mathematics and Mechanics, Kyiv Taras
Shevchenko University, 64 Volodymyrska st., 01033 Kyiv, UKRAINE,
\vspace{0.3cm}

\noindent
V.Maz.: Department of Mathematics, Uppsala University, Box. 480,
SE-75106, Uppsala, SWEDEN, email: {\tt mazor\symbol{64}math.uu.se}

\end{document}